\documentclass[12pt]{article}
\usepackage{amsfonts,amssymb,latexsym,
amscd,
theorem,epsfig,graphics}
\newtheorem{lemma}{Lemma}[section]
\newtheorem{proposition}[lemma]{Proposition}
\newtheorem{remark}[lemma]{Remark}
\newtheorem{example}[lemma]{Example}
\newtheorem{theorem}{Theorem}
\newtheorem{definition}[lemma]{Definition}
\newtheorem{corollary}[theorem]{Corollary}
\begin{document}
\newcommand{\eps}{{\varepsilon}}
\newcommand{\proofend}{$\Box$\bigskip}
\newcommand{\sign}{{\mbox{\rm sign}}}
\newcommand{\mes}{{\mbox{\rm mes}}}
\newcommand{\lcm}{{\mbox{\rm lcm}}}
\newcommand{\Id}{{\mbox{\rm Id}}}
\newcommand{\tet}{{\theta}}
\newcommand{\C}{{\mathbf C}}
\newcommand{\Q}{{\mathbf Q}}
\newcommand{\R}{{\mathbf R}}
\newcommand{\Z}{{\mathbf Z}}
\newcommand{\cc}{{\cal C}}
\renewcommand{\k}{\mathbf k}
\newcommand{\coker}{{\rm Coker}}
\newcommand{\RP}{{\mathbf {RP}}}
\newcommand{\rk}{{\rm rk}}
\newcommand{\kk}{{\mathbf k}}

\title
{Topology of cyclic configuration spaces and
periodic trajectories of multi-dimensional billiards}
\author{Michael Farber
\thanks{Partially supported by a grant from the
Israel Academy of Sciences and Humanities and by
the Herman Minkowski Center for Geometry
\newline
${}\quad ^{**}$Partially supported by an NSF grant.}
$\ $ and
Serge Tabachnikov$^{**}$ \\
{\it School of Mathematical Sciences, Tel Aviv University}\\
{\it Ramat Aviv, Tel Aviv 69978, Israel}\\
and\\
{\it Department of Mathematics, University of Arkansas}\\
{\it Fayetteville, AR 72701, USA}\\
e-mail: {\it farber@ math.tau.ac.il, serge@comp.uark.edu}}
\date{November 25, 1999}
\maketitle
\begin{abstract} We give lower bounds on the number of periodic trajectories in
 strictly convex smooth billiards in $\R^{m+1}$ for $m\geq 3$.
For plane billiards (when $m=1$) such bounds were obtained by G. Birkhoff
in the 1920's.
Our proof is based on topological methods of calculus of variations --
equivariant Morse and
Lusternik-Schnirelman theories. We compute the equivariant cohomology ring
 of the cyclic
configuration space of the sphere $S^m$, i.e., the space of $n$-tuples of
points  $(x_1,\dots, x_n)$, where $x_i\in S^m$ and $x_i \neq x_{i+1}$ for
$i=1,\dots,n$.
\smallskip

\end{abstract}

\section{Introduction}

The billiard dynamical system describes the free motion of a mass-point in
a domain in Euclidean space with a reflecting boundary: a point moves along
a straight line with unit speed until it hits the boundary, at the impact
point the normal component of the velocity instantaneously changes sign
while the tangential component remains the same, and the rectilinear motion
continues with a new unit velocity. We refer to \cite {KT}, \cite {Ta} for
surveys of mathematical billiards.

We address the following problem:
how many periodic billiard trajectories are there in a smooth strictly
convex domain in ${\R}^{m+1}$?

For plane billiards ($m=1$)
this problem was studied by G. Birkhoff in \cite {Bi}. A closed billiard
orbit of a plane billiard is an inscribed plane
$n$-gon whose consecutive sides make equal angles with a closed convex
plane curve $X$.
Such a periodic orbit has a rotation number
$0<r<n/2$. Birkhoff proved that {\it for every $n \geq 2$ and $r \leq n/2$,
coprime
with $n$, there
exist at least two distinct $n$-periodic billiard trajectories with the
rotation number $r$} -- see \cite {Bi}. From contemporary view point, this
result
follows from the theory of area-preserving twist maps -- see, e.g.,
\cite {Ban} for a survey.

In the present paper we use a well-known
variational reduction of the periodic billiard trajectories problem
which is based on the observation that these trajectories are critical
points of the perimeter length functional on the variety of inscribed polygons.
This allows to apply topological methods (Morse and
Lusternik-Schnirelman theories)
and reduces the problem to studying the topology
of {\it the cyclic configuration space}. Inscribed $n$-gons
are in one-to-one correspondence with sequences
$$(x_1, x_2, \dots, x_n)\in X^{\times n}= X \times \dots \times X,$$ such that
\[x_{1}\ne x_{2}, \, x_{2}\ne x_{3}, \, \dots, x_{n}\ne x_{1};\]
here $X\subset \R^{m+1}$ denotes the boundary of the billiard domain,
topologically, the sphere.
The variety of all such configurations is the {\it cyclic configuration
space}. Our main effort in this paper is in
studying the cohomology ring of this configuration space. Note
that the length functional, considered as a map $X^{\times n}\to \R$ on the
total Cartesian
power of $X$, fails to be smooth at the points with $x_i=x_{i+1}$ for some $i$.
This explains why the cyclic configuration space is the natural
topological object related to the billiard problem.

The dihedral group $D_n$ acts on the cyclic configuration space; its action
is generated by the cyclic permutation and the reflection
\[(x_1, x_2, \dots, x_n)\mapsto (x_2, x_3, \dots, x_n, x_1),\quad
(x_1, x_2, \dots, x_n)\mapsto (x_n, x_{n-1}, \dots, x_1).\]
Two $n$-periodic billiard trajectories will be considered the
same if they belong to the same orbit of $D_n$. In particular, the estimates
in Theorem \ref {thm1} below concern the number of distinct $D_n$-orbits in
the cyclic configuration space.

Our main result is as follows.

\begin{theorem} \label {thm1}
Let $X\subset \R^{m+1}$ be a smooth
strictly convex hypersurface, where $m \geq 3$. Fix an odd number $n\ge 3$.
Then

(A) The number of distinct $D_n$-orbits of  $n$-periodic billiard trajectories
inside $X$ is not less than
\begin{eqnarray}
[\log_2(n-1)]+m.\label{formula1}
\end{eqnarray}

(B) For a generic $X$, the number of distinct $D_n$-orbits of
$n$-periodic billiard trajectories  inside $X$ is not less than
\begin{eqnarray}
(n-1)m\label{formula2}.
\end{eqnarray}
\end{theorem}

Square  brackets in  (\ref{formula1}) denote the floor function, i,e., the
largest integer  not
exceeding $\log_2(n-1)$.

We deduce statement (A) of Theorem \ref{thm1} from
Lusternik-Schnirelman theory and statement (B) from Morse theory; this explains
the genericity assumption in (B) (precise definition of the needed
genericity is given in Section 4).

We conjecture that the estimate in case (A) can be
improved to be
$n+m-2$.

It is important to note that if the period $n$ in Theorem \ref{thm1} is not
prime then
$n$-periodic billiard trajectories, whose existence is asserted by this
theorem, may be
multiple ones; for example, instead of a genuine $9$-periodic orbit one
might have a $3$-periodic
orbit, traversed three times.

To the best of our knowledge, the only previous attempt to extend
Birkhoff's results to
multi-dimensional setting was made by I. Babenko in \cite {Ba}. He applied
the variational
approach to periodic trajectories of convex billiards in 3-dimensional space.
Unfortunately, the paper \cite {Ba} contains an error (cf. Remark \ref{rmk34}), and the situation
with $m=2$ remains unclear.
We would like to emphasize that ideas and results from \cite{Ba} are
substantially used
in the present paper.

We would like to mention, in passing, the case of 2-periodic billiard
trajectories. These are the
diameters of the billiard hypersurface $X \subset \R^{m+1}$, that is, the
chords, perpendicular to
the hypersurface at both ends. It is well known that the least number of
diameters is $m+1$, and this
estimate is sharp as the example of a generic ellipsoid shows.

The arguments that prove Theorem ~\ref{thm1} apply to a much wider class of
billiards. Namely, let
$X$ be  a smooth strictly geodesically convex billiard hypersurface in a
Riemannian manifold such
that there exists a unique geodesic line through every two points of $X$.
For example, this
condition holds when the  ambient manifold is the spherical or the
hyperbolic space or, more
generally, a Hadamard space.

\begin{theorem} \label {thm1.5}
For $X$ as in the preceding paragraph, both statements of Theorem
~\ref{thm1} hold true, provided $\dim X \ge 3$.
\end{theorem}

The content of the paper is, briefly, as follows. In Section 2 we construct
a general spectral sequence that
computes the cohomology ring  of the cyclic configuration space of a smooth
manifold.  In  Section 3 we use it to compute the cohomology ring  of the
cyclic
configuration space of the sphere. Section 4 is Morse theoretical:
we describe the relations (mostly known) between the information about the
periodic billiard trajectories
and the topology of the cyclic configuration space
of the billiard hypersurface. In the last Section 5 we compute the
$D_n$-equivariant cohomology  of
the cyclic configuration space of the sphere; the computation makes use of
the topological
results from Section 3 and the Morse-theoretical ones from Section 4. We
finish the paper with a
proof of our main result, Theorems ~\ref{thm1} and ~\ref{thm1.5}.
\medskip

{\bf Acknowledgments}. It is a pleasure to acknowledge the hospitality of
Max-Planck-Institut in
Bonn where the authors started to work on this paper. Many thanks to Burt
Totaro for patiently explaining to us his paper \cite {To} and to Ivan
Babenko for discussions on his paper \cite{Ba}.

\section{Cyclic configuration space of a smooth manifold}

 Let $X$ be a smooth manifold. Consider the
configurations of ordered
 $n$-tuples $(x_1, x_2, \dots, x_n)$ of points of $X$ satisfying:
\begin{eqnarray}
x_i\ne x_{i+1}\quad {\rm for}\quad i=1, 2, \dots, n.\label{cyclic}
\end{eqnarray}
Here and elsewhere below we understand the indices cyclically, that is
$n+i= i$ (in particular,
(\ref{cyclic}) contains the requirement $x_n\ne x_1$).
The space of all such configurations will be called the {\it cyclic
configuration space} and denoted
by $G(X,n)$. The dihedral group $D_n$ acts naturally on $G(X,n)$ (cf. above);
this action is free if $n$ is an odd prime.

\begin{example} \label{ex1}
{\rm The simplest example of a cyclic configuration space is provided by
$G(S^1,n)$,
the cyclic configuration space of a circle, which plays an important role
in Birkhoff's theory of
convex plane billiards \cite{Bi}.
A configuration $(x_1, x_2, \dots, x_n)\in G(S^1,n)$ can be uniquely
described by
the initial point $x_1\in S^1$ and by the angles $0<\phi_i<1$, where $i=1,
2, \dots, n-1$,
assuming that the sum
$\phi_1+\phi_2+\dots +\phi_{n-1}$
is not an integer. Namely, given $x_1$ and $\phi_1, \dots, \phi_{n-1}$ we set
\[x_2=x_1\exp(2\pi i\phi_1),\, \,
 x_3=x_2\exp(2\pi i\phi_2),\, \,
\dots, x_n=x_{n-1}\exp(2\pi i\phi_{n-1}).\]
Here we identify $S^1$ with the unit circle in $\C$. The hyperplanes
\[\phi_1+\phi_2+\dots +\phi_{n-1} = r, \quad {\rm where}\quad  r=1, 2,
\dots, n-2,\]
divide the cube $(0,1)^{n-1}$ into domains $N_1, N_2, \dots, N_{n-1}$, where
\[N_r=\{(\phi_1, \phi_2, \dots, \phi_{n-1}):  r-1< \sum_{i=1}^{n-1}\phi_i <
r\}.\]
Each $N_r$ is homeomorphic to $\R^{n-1}$.
It follows that the cyclic configuration space $G(S^1, n)$ is homeomorphic
to a disjoint union
of $S^1\times N_r$, where $r=1, 2,\dots, n-1$. Note that the number $r$,
determining the
connected component of a given configuration, is precisely the rotation
number. Let us also indicate
how the dihedral group $D_n$ acts on $G(S^1,n)$. Given a point
$(x_1, \phi_1, \dots, \phi_{n-1})\in S^1\times N_r$,
the cyclic permutation takes it to
\[(x_1\exp(2\pi i\phi_1), \phi_2, \dots, \phi_{n-1},
r-\sum_{i=1}^{n-1}\phi_i).\]
The reflection maps it to
\[(x_1, 1-r +\sum_{i=1}^{n-1}\phi_i, 1-\phi_{n-1}, 1-\phi_{n-2}, \dots,
1-\phi_2).\]
Hence the cyclic permutation preserves the rotation number, and
the reflection maps configurations with
rotation number $r$ to configurations with rotation number $n-r$.}
\end{example}

The cyclic configuration space $G(X,n)$
is not to be confused with the usual  configuration space
$$F(X,n) = \{ (x_1, \dots , x_n) \in X^{\times n}: x_i \neq x_j \quad {\rm
for} \quad i \neq j \}.$$
Recall a description of the cohomology ring $H^*(F({\R}^m,n);{\Z})$,
obtained by F. Cohen -- see
\cite {Co1}, \cite {Co2}.
The space $F({\R}^m,2)$ is homotopy equivalent to $m-1$-dimensional sphere;
denote by $\omega$
its top-dimensional cohomology class. Let $p_{ij}$
be the projection  $F({\R}^m,n) \to F({\R}^m,2)$ on the $i$-th and $j$-th
components, and let
$G_{i,j} = p_{ij}^* (\omega)$. Then the ring $H^*(F({\R}^m,n);{\Z})$, where
$m>1$,
is the graded-commutative algebra over ${\Z}$ with generators
$$G_{i,j}, \ 1 \leq i,j \leq n, \ i \neq j, \quad {\rm deg}\ G_{i,j} = m-1$$
and relations
\begin{enumerate}
\item[(a)] $G_{i,j} = (-1)^m G_{j,i}$,
\item[(b)] $G_{i,j}^2 =0$,
\item[(c)] $G_{i,j}G_{i,k} + G_{j,k}G_{j,i} + G_{k,i}G_{k,j} = 0$ \ for
$i,j,k$ distinct.
\end{enumerate}
An additive basis for $H^{r(m-1)}(F({\R}^m,n);{\Z})$ is given by monomials
$G_{i_1, j_1}G_{i_2, j_2} \dots G_{i_r, j_r}$ with
$i_1 < \dots < i_r$ and $i_k < j_k$ for all $k$.

We start the study of topology of cyclic configuration spaces with the
space $G({\R}^m,n)$.
The following proposition is an analog of Cohen's result.

\begin{proposition} \label{prop1} For $m>1$
the ring $H^*(G({\R}^m,n);\Z)$ is a graded commutative algebra over ${\Z}_2$
with  generators
$$s_1, s_2, \dots s_n, \quad {\it with}\quad {\rm deg}\ s_i = m-1$$
and the relations
\[s_1^2=s_2^2=\dots =s_n^2=0,\]
\begin{eqnarray}
s_1s_2\dots s_{n-1}+\epsilon s_2s_3\dots s_n + \epsilon^2 s_3s_4\dots
s_ns_1+\dots +\epsilon^{n-1}
s_ns_1\dots s_{n-2}=0,
\label{relation}
\end{eqnarray}
where $\epsilon = (-1)^{(m-1)(n-1)}$.
If $\phi$ denotes the inclusion of $F({\R}^m,n)$ into $G({\R}^m,n)$, then
$\phi^* (s_i) = G_{i, i+1}$.
\end{proposition}

{\bf Proof.}
In general, the cohomology groups of the complement ${\R}^N - \cup V_i$
of a collection of affine subspaces $V_i\subset \R^N$, where $i=1,\dots, k$,
is described in \cite {GM}. Recall this description, following \cite{Va},
chapter 3, \S 6. For $I \subset \{1,\dots,k\}$ let
$$V_I = \cap_{i \in I} V_i.$$
Two subsets $I, J\subset \{1,\dots,k\}$ are equivalent, $I \sim J$,  if
$V_I =  V_J$.
Each equivalence class has a unique  maximal element. Given such a maximal
element $I$,
consider the standard chain complex of the simplex $\Delta_I$ spanned by $I$.
Factorize this complex by the subcomplex generated by the faces
$\Delta_J$ with $J \subset I, J \nsim I$, and denote the resulting quotient
complex by $C(I)$.
Then the cohomology of the complement ${\R}^N - \cup V_i$ is given by
$$\tilde H^q ({\R}^N - \cup V_i;\Z) = \oplus \tilde H_{N-{\rm dim}\ V_I -q
-1}(C(I);\Z),$$
the sum taken over all the equivalence classes with maximal elements $I$,
such that $V_I$ are nonempty.

In our situation, ${\R}^N = ({\R}^{m})^{\times n}$, and $V_i = \{x_i =
x_{i+1}\}$ for
$i=1, 2, \dots, n$. If $|I| \leq n-2$ then the equivalence class of $I$
consists  of one element only;
on the other hand, all subsets of cardinality $n-1$ are equivalent to each
other and to $\{1,\dots,n\}$.
Therefore, if $|I| \leq n-2$ then
$$\tilde H_{\ast}(C(I);\Z) = H_{\ast}(\Delta_I, \partial \Delta_I;\Z) =
\tilde H_{\ast} (S^{|I|-1};\Z),$$
and each subset $I$ with $s=|I|  \leq n-2$ makes a contribution of a copy
of $\Z$ to
$H^{s(m-1)}(G({\R}^m,n);\Z)$. Similarly, if $I = \{1,\dots,n\}$ then
$$\tilde H_{\ast}(C(I);\Z) = H_{\ast}(\Delta_I, sk_{n-2} (\Delta_I);\Z) =
 \tilde H_{\ast} (\vee_{i=1}^{n-1} S_i^{n-1};\Z),$$
(where $sk_{n-2} (\Delta_I)$ denotes the $(n-2)$-dimensional skeleton of
the simplex $\Delta_I$)
and $I$ makes a contribution of $\Z^{n-1}$ to
$H^{(n-1)(m-1)}(G({\R}^m,n);\Z)$.
This implies that $\dim H^{\ast}(G({\R}^m,n);\Z)$  is free abelian of rank
\[\rk\ H^{s(m-1)}(G({\R}^m,n);\Z) = \left\{
\begin{array} {r@{\quad \quad }l}
\displaystyle{n \choose s}, & {\rm for} \quad 0\le s\le n-2,\\ \\
n-1, &  {\rm for} \quad s=n-1
\end{array}
\right.
\]
Hence, the Poincar\'e polynomial of the cyclic configuration
space $G({\R}^m,n)$  equals
$$(t^{m-1} + 1)^n - t^{(n-1)(m-1)} - t^{n(m-1)}.$$

For $i=1, \dots, n$ denote by $s_i$ the generators of $H^{m-1}(G({\R}^m,n);\Z)$
coming from 1-element sets $\{i\}$.
It is clear that $s_i$ is the pull-backs of the top-dimensional class of
$G({\R}^m,2) = F({\R}^m,2) \approx S^{m-1}$ under the projection
$p_{i,i+1}$. Therefore $s_i^2=0$
and $\phi^* (s_i) = G_{i, i+1}$.

It remains to show that an additive basis for $H^{r(m-1)}(G({\R}^m,n);\Z)$
with $r \leq n-1$,
is given by the monomials
$s_{i_1} \dots s_{i_r}$, where $i_1 < \dots < i_r$, with one additional
relation (\ref{relation}).

Denote by $S=S^{m-1}\subset \R^m$ the unit sphere and let $S^{\times
(n-1)}$ be the
Cartesian power of $S$. For any $p=1, 2, \dots, n$ we construct an embedding
\[g_p: S^{\times (n-1)}\to G(\R^m,n),\quad (y_1,y_2, \dots, y_{n-1})\mapsto
(x_1, x_2, \dots, x_n),\]
where
$$
\begin{array}{l}
x_p=0,\\
x_{p+1}=y_1,\\
x_{p+2}=y_1+y_2,\\
\dots\\
x_n = y_1+y_2+\dots + y_{n-p}\\
x_1=y_1+y_2+\dots +y_{n-p}+y_{n-p+1}\\
\dots\\
x_{p-2}=y_1+y_2+\dots +y_{n-2}\\
x_{p-1}=y_1+y_2+\dots +y_{n-2}+ny_{n-1}.\end{array}
$$
If $\overline s_i\in H^{m-1}(S^{\times (n-1)};\Z)$ denotes the obvious
generator
corresponding to the $i$-th factor (where $i=1, 2, \dots, n-1$),
we have
$$
g_p^\ast(s_i) = \left\{
\begin{array}{l}
\overline s_{n-p+1+i},\quad{\rm for}\quad 1\le i\le p-2,\\
\overline s_{i-p+1}, \quad{\rm for}\quad p\le i\le n,\\
 (-1)^m \overline s_{n-1},\quad {\rm for}\quad i=p-1.
\end{array}
\right.
$$
All but the last relation being obvious, let us explain the relation
$g_p^\ast(s_{p-1})=(-1)^m \overline s_{n-1}$. Consider the homotopy
$F_t: S^{\times (n-1)}\to S$, where $t\in [0,1],$
given by
\[F_t(y_1, \dots, y_{n-1}) = \frac{ty_1+\dots +ty_{n-2}+ny_{n-1}}
{||ty_1+\dots +ty_{n-2}+ny_{n-1}||},\quad y_i\in S.\]
It is clear from the definitions that $F_1^\ast(v)=(-1)^m
g_p^\ast(s_{p-1})$, where
$v\in H^{m-1}(S;\Z)$ is the
generator. On the other hand, the map $F_0$ is just the projection
$(y_1, \dots, y_{n-1})\mapsto y_{n-1}$ and hence $F_0^\ast(v)=\overline
s_{n-1}$.

Let $a_i\in H^{(n-1)(m-1)}(G(\R^m,n);\Z)$ denote the product $s_{i+1}
s_{i+2} \dots s_n s_1
\dots s_{i-1}$. Then we obtain
$$
g_p^\ast(a_i) =\left\{
\begin{array}{l}
\overline s_1\overline s_2\dots \overline s_{n-1},\quad{\rm for}\quad i =p-1,\\
- \epsilon\overline s_1\overline s_2\dots \overline s_{n-1},\quad{\rm
for}\quad i =p-2,\\
0, \quad {\rm for} \quad i\ne p-2, \, p-1,
\end{array}
\right.
$$
where $\epsilon = (-1)^{(m-1)(n-1)}$.
Varying $p$, this shows that all $a_i$ are nonzero and also that for
any nontrivial relation $\sum \beta_i a_i$
in $H^{(n-1)(m-1)}(G(\R^m,n);\Z)$ one must have $\beta_{i+1} =\epsilon
\beta_i$.
In other words, any relation of degree $n-1$ between the classes $s_i$ must
be a consequence of
(\ref{relation}). The rank calculation above shows that at least one
nontrivial relation between
the classes
$a_1, \dots, a_n$ exists, and so (\ref{relation}) holds.

Let us finally show that there exist no nontrivial relations between $s_1,
\dots, s_n$ of degree
$r<n-1$.  Assume that there exists such a relation $x$ containing a
monomial $s_{i_1}s_{i_2}\dots s_{i_r}$. Choose a set of indices $\{j_1,
\dots, j_{n-r-1}\}$,
disjoint from the set $\{i_1,\dots,i_r\}$, and let $y=s_{j_1}\dots
s_{j_{n-r-1}}$. Then
$xy$ is a relation of degree
$n-1$, containing neither of the terms $a_{j_1},\dots, a_{j_n-r-1}$;  we
have shown above that it is
impossible. This completes the proof. {\proofend}

\begin{remark}\label{rmk1} {\rm The cyclic configuration space
$G(\R^m,3)$ coincides with $F(\R^m,3)$. In this case relation (\ref{relation})
turns into F. Cohen's
relation (c) in the cohomology $H^\ast(F(\R^m,3);\Z)$ after substituting
$s_i=G_{i,i+1}$.}
\end{remark}

Next, we construct a spectral sequence computing the cohomology  of the
cyclic configuration space
$G(X,n)$, where $X$ is an arbitrary smooth orientable manifold. This
spectral sequence
is an analog of the one constructed by Totaro in \cite {To}
for the configuration space $F(X,n)$.

Denote by
$p_j: X^{\times n} \to X$ the projection on the $j$-th component and by
 $q_j=p_{j, j+1}$ the projection $ X^{\times n} \to X^{\times 2}$ on the
$j$-th and $j+1$-th
components. Also, let $\Delta \in H^m (X^{\times 2};\kk)$ denote the
cohomology class of the
diagonal; here $m=\dim X$ and $\kk$ is a field.

\begin{theorem}  \label {thm2} Let $X$ be a smooth orientable manifold of
dimension $m>1$;
let $\kk$ be a field.

(A) There exists a spectral sequence of bigraded differential algebras over
$\kk$ which converges to $H^\ast(G(X,n);\kk)$ whose $E_2$-term is the
quotient of the bigraded commutative algebra
\[H^\ast(X^{\times n};\kk) \otimes H^\ast(G({\R}^m,n);\kk),\]
where $H^p(X^{\times n};\kk)$ has
bidegree $(p,0)$ and $H^q(G(\R^m,n);\kk)$ has bidegree $(0,q)$, by the
relations
 $$p_i^\ast(v)s_i = p_{i+1}^\ast(v)s_i,\ i=1, \dots, n,$$
where $v\in H^\ast(X;\kk)$ is an arbitrary class.

(B) The action of the dihedral group $D_n$  on the spectral
sequence is given by the action on $H^\ast(X^{\times n};\kk)$,
induced by its action on $X^{\times n}$, and by $\tau(s_i) = s_{\tau(i)}$.

(C) The first non-trivial differential is $d_m$, where $m=\dim X$. It is
defined by the formulas
$$d_m s_i = q_i^* (\Delta) \quad {\rm and} \quad d_m (H^\ast(X^{\times
n};\kk)) = 0.$$
\end{theorem}

{\bf Proof.} The proof is a modification of the arguments  by B. Totaro in
\cite {To}.

Consider the inclusion $\psi: G(X,n) \to X^{\times n}$ and the Leray
spectral sequence of the continuous map $\psi$:
$$E_2^{p,q}\, =\, H^p(X^{\times n}; R^q \psi_{\ast} {\kk}) \Rightarrow
H^{p+q} (G(X,n);\kk)$$
where $R^q \psi_{\ast} {\kk}$ is the sheaf on $X^{\times n}$ associated
with the presheaf
$$U \mapsto H^q(U \cap G(X,n);\kk).$$

We will consider partitions of the set of indices $\{1,\dots,n\}$ into
intervals, that is,
subsets of the form $\{i,i+1,i+2, \dots, i+j\}$; as usual, the indices are
understood cyclically.
For example, the following are interval partitions of the set $\{1, 2, 3,
4, 5\}$:
\[\{1,2\}\cup \{3, 4, 5\}, \quad \{2\}\cup\{3\}\cup\{4, 5, 1\}.\]
Let $J$ be such a partition; denote by $X_J$ the subset of $X^{\times n}$
consisting
of configurations $c=(x_1, x_2, \dots, x_n)\in X^{\times n}$
satisfying the conditions $x_i = x_j$ if $i$ and $j$ lie in the same
interval of the partition $J$.
Given two interval partitions $I$ and $J$, we  say that $J$ {\it
refines} $I$ and write $I\prec J$
if the intervals of $I$ are unions of the intervals of $J$.
Denote by $|J|$ the number of intervals in the partition $J$.
$X_J$ is naturally identified with the Cartesian power $X^{\times |J|}$.
The case $|J|=1$ corresponds to the deepest diagonal
$\{(x,x,\dots, x)\}\subset X^{\times n}$.
Note that $I\prec J$ implies $X_I\subset X_J$ and $|I|\le |J|$.

Denote by $D(X,n)$  the subset of $X^{\times n}$ satisfying the conditions
$x_i \neq x_{i+1}$ for $i=1, \dots, n-1$, but not necessarily for $i=n$.
It is clear that $D({\R}^m,n)$ is homotopy equivalent to the product
$(S^{m-1})^{\times (n-1)}$
and so its cohomology admits a description similar to Proposition \ref{prop1},
but without relation (\ref{relation}).

Let $J$ be a partition of $\{1, 2, \dots, n\}$ on intervals of lengths
$j_1,...,j_r$
as above, and let $c$ be a configuration with
$c \in X_J$ such that $c$ does not belong to any $X_I$ if
$J$ refines $I$.
We claim that the stalk of the sheaf $R^q \psi_{\ast} {\kk}$ at $c$ equals
$$
(R^q \psi_{\ast} {\kk})_c =
\left\{
\begin{array}{l}
H^q(D({\R}^m,j_1) \times \dots \times D({\R}^m,j_r);\kk) \quad {\rm if} \quad r
\geq 2\\ \\
H^q(G({\R}^m,n);\kk) \quad {\rm if} \quad r=1.
\end{array}
\right.
$$
Indeed, by definition, this stalk is $H^q(U \cap G(X,n);\kk),$
 where $U$ is a small open ball around $c$. If $c=(x_1, x_2, \dots x_n)$
then we may choose points
$y_1, y_2, \dots, y_r\in X$, one for each interval of $J$, so that
$x_i=y_{j_s}$ if $i$ belongs to the
$j_s$-th interval.  Let $U_j\subset X$ be a small open neighborhood of $y_j$,
so that each $U_j$ is diffeomorphic to $\R^m$ and
the sets $U_j$ and $U_{j'}$ are disjoint when the points $y_j$
and $y_{j'}$ are distinct. Then we may take
$U=U_1^{\times j_1}\times U_2^{\times j_2}\times \dots \times U_r^{\times
j_r}$,
and our claim follows.

In particular, we see that $R^q \psi_{\ast} {\kk}$ vanishes unless $q$ is
a multiple of $m-1$.

If $|J|\geq 2$ then
$$(R^* \psi_{\ast} {\kk})_c = H^* ((S^{m-1})^{\times (n-|J|)};\kk),$$
and hence
\[{\rm dim}\ (R^{s(m-1)} \psi_{\ast} {\kk})_c = \displaystyle{n-|J|
\choose s}.\]
Note that there is a unique top-dimensional class of dimension
$(m-1)(n-|J|)$.

If $|J|=1$ (which means that $X_J$ is the deepest diagonal) then, by
Proposition~\ref{prop1},
$${\rm dim}\ (R^{r(m-1)} \psi_{\ast} {\kk})_c =
\left\{
\begin{array}{l}
\displaystyle{n \choose r} \quad {\rm for} \quad r \leq n-2,\\ \\
 n-1\quad {\rm for}\quad r=n-1,
\end{array}
\right.
$$
and hence there exists $(n-1)$-dimensional space of top-dimensional
cohomology classes; they have dimension $(n-1)(m-1)$.
Note that the sheaf $R^{(n-1)(m-1)}\psi_\ast\kk$ is a constant sheaf
$\kk^{n-1}$
with support on the deepest diagonal $X_J$. This follows repeating the
arguments of Totaro \cite{To}, p. 1062. The fact that this sheaf is locally
constant is easy
since this sheaf can be represented as sheaf of cohomology of fibers of a
fibration over $X_J$.
A further statement that this sheaf is constant (i.e., has a trivial monodromy)
follows from explicit calculation of the cohomology $H^\ast(G(\R^m,n);\kk)$
(which serves as the cohomology of the fiber),
since we may explicitly label the cohomology classes by
polynomials in $s_1, \dots, s_n$.

We wish to have a similar description for the sheaves
$R^{r(m-1)}\psi_\ast\kk$ with $r<n-1$.
We will show that they can be represented as direct sums of constant
sheaves supported on
different diagonals $X_J$.
For an interval partition $J$ of $\{1, 2,\dots, n\}$ with $|J|>1$
 denote by $\varepsilon_J$ the constant sheaf
with stalk $\kk$ and support $X_J$. We claim that {\it for $r<n-1$ the sheaf
$R^{r(m-1)} \psi_{\ast} {\kk}$ is isomorphic to the direct sum of sheaves
\begin{eqnarray}
R^{r(m-1)} \psi_{\ast} {\kk}\, \simeq\, \bigoplus_{|J|=n-r}
\varepsilon_J,\label{split}\end{eqnarray}the sum taken over all interval
partitions $J$ with $
|J|=n-r$.}

To prove the above claim we first make the following remark. Let
$I$ be an interval partition
of $\{1, 2, \dots, n\}$ into intervals of length $i_1, i_2, \dots, i_s$,
where $s=|I|>1$. Then for any
interval partition $J$ into intervals of length $j_1, j_2, \dots, j_{n-r}$,
such that $I\prec J$, we have the
canonical inclusion
\[\nu_{JI}: D(\R^m,i_1)\times \dots \times D(\R^m,i_s) \to
D(\R^m,j_1)\times \dots \times D(\R^m,j_{n-r}).\]
Note that the target of $\nu_{JI}$ has a unique nonzero
$r(m-1)$-dimensional class.
The induced map on $r(m-1)$-dimensional cohomology with $\kk$ coefficients
\[\nu^\ast_{JI}: H^{r(m-1)}(D(\R^m,j_1)\times \dots\times
D(\R^m,j_{n-r}))\to
H^{r(m-1)}(D(\R^m,i_1)\times \dots\times D(\R^m,i_{s}))\]
is a monomorphism. Let $z_{JI}$ denote the image of the top-dimensional
generator.
Then we observe (similarly to Lemma 3 in \cite{To})
that for a fixed $I$ the classes $\{z_{JI}\}$ form a linear basis of the
cohomology
$H^{r(m-1)}(D(\R^m,i_1)\times \dots\times D(\R^m,i_{s});\kk)$, where $J$ runs
over all partitions with
$I\prec J$  and $|J|=n-r$.

A similar statement holds in the case $|I|=1$.
Namely, for any interval partition $J$ into intervals of length $j_1, j_2,
\dots, j_{n-r}$,
where $n-r>1$, we have the canonical inclusion
\[\nu_{JI}: G(\R^m,n) \to D(\R^m,j_1)\times \dots \times D(\R^m,j_{n-r});\]
 denote by $z_{JI}$ the image of the top-dimensional generator under
the induced map on the cohomology,
$z_{JI}\in H^{r(m-1)}(G(\R^m,n);\kk)$. Then the set $\{z_{JI}\}$ (where $J$
runs over all partitions
with $|J|=n-r$) forms a linear basis of
$H^{r(m-1)}(G(\R^m,n);\kk)$. This follows from the explicit calculation of the
cohomology of $G(\R^m,n)$,
given in Proposition \ref{prop1}.

Given an interval partition $J$ on intervals of length $j_1, \dots,
j_{n-r}$, where $n-r>1$,
consider the following commutative diagram
$$
\begin{array}{clc}
G(X,n) & \to & X^{\times n}\\
\downarrow & & \downarrow {\rm id}\\
D(X,j_1)\times \dots \times D(X,j_{n-r}) & \stackrel {g_J}{\to} & X^{\times n}
\end{array}
$$
formed by the natural inclusions. Define a sheaf $\varepsilon'_J$ over
$X^{\times n}$ by
$$\varepsilon'_J = R^{r(m-1)}{g_J}_\ast(\kk).$$
We want to show that $\varepsilon'_J$ is isomorphic to $\varepsilon_J$,
defined above. First,
$\varepsilon'_J$
vanishes outside $X_J$ (since we are considering the cohomology of the top
dimension). Let
$U$ be a small open neighborhood of a point $c\in X_J\subset X^{\times
n}$, such that
$U=\prod U_i$, where all $U_i$ are small
open disks and $U_i=U_j$ if $i$ and $j$ lie in the same interval
of $J$. Then
$$\varepsilon'_J(U) = H^{r(m-1)}(D(U_{i_1}, j_1)\times \dots \times
D(U_{i_{n-r}},j_{n-r});\kk) =\kk.$$
Hence  $\varepsilon'_J$ is a constant sheaf with stalk $\kk$
supported on $X_J$, i.e.,
$\varepsilon'_J\simeq\varepsilon_J$.

The commutative diagram above gives a map of sheaves $\varepsilon_J\to
R^{r(m-1)}\psi_\ast(\kk)$
and  we obtain a map of sheaves
\[\bigoplus_{|J|=n-r} \varepsilon_J\, \to\,  R^{r(m-1)}\psi_\ast(\kk),\]
which, as we have seen above, is an isomorphism on stalks;
hence it is an isomorphism, and the claim (\ref{split}) follows.

We arrive at the following description  of the term $E_2$ of the spectral
sequence:
$$
E_2^{p, r(m-1)} =
\left\{
\begin{array}{l}
H^p(X_{J_0};\kk) \otimes \kk^{n-1},\quad{\rm for}\quad r=n-1,\quad{\rm
where}\quad |J_0|=1,\\ \\
\bigoplus_{|I|=n-r} H^p(X_I;\kk),\quad{\rm for}\quad r<n-1
\end{array}
\right.
$$
We will now identify this description with the one given in the statement
of the theorem.
Consider first the case $r < n-1$.
Assign to a monomial $s_{i_1} \dots s_{i_r}$ with $i_1< i_2<\dots <i_r$
the equivalence relation on the set of indices $\{1,\dots,n\}$ generated by
the relations:
$$i_1 \sim i_1 +1, \, \, i_2 \sim i_2 + 1, \dots, i_r \sim i_r + 1.$$
This equivalence relation defines a partition $I$ of the set
$\{1,\dots,n\}$ on $n-r$ intervals.
In view of the relations $p_i^\ast(v)s_i = p_{i+1}^\ast(v)s_i$, the term
$H^p (X^{\times n};\kk) s_{i_1} \dots s_{i_r}$ is isomorphic to $H^p
(X_I;\kk)$,
and we are done in this case.
Consider now the case $r=n-1$. There are $n-1$ linearly
independent monomials of degree $n-1$ in classes $s_1, \dots, s_n$,
and for every such monomial the corresponding interval partition equals
$J_0 = \{1,\dots,n\}$.
Hence the image of
$H^p (X^{\times n};\kk) \otimes H^{(n-1)(m-1)} (G({\R}^m,n);\kk)$
is isomorphic to $H^p (X_{J_0};\kk) \otimes {\kk}^{n-1}$ after imposing the
relations
$p_i^\ast(v)s_i = p_{i+1}^\ast(v)s_i$, as stated.

Now we consider the differentials of the spectral sequence. The first $m-2$
of them, $d_2, \dots, d_{m-1}$, vanish by the dimension considerations. To
find $d_m$ consider the inclusion
$\phi: F(X,n) \hookrightarrow G(X,n)$. This inclusion induces a
homomorphism of the spectral sequence for $G(X,n)$ to that for $F(X,n)$,
constructed in \cite {To}. In the latter spectral sequence the $E_2$-term
is the quotient of the graded commutative algebra
$$H^\ast(X^{\times n};\kk) \otimes H^\ast(F({\R}^m,n);\kk)$$
modulo the relations
 $$p_i^\ast(v) G_{ij} = p_j^\ast(v) G_{ij}\quad {\rm for} \quad i\neq
j\quad {\rm and} \quad v\in H^\ast(X;\kk),$$
and the first non-trivial differential acts as follows:
$d_m G_{ij} = p^*_{ij} \Delta.$
According to Proposition~\ref{prop1}, $\phi^*(s_i) = G_{i,i+1}$, and,
obviously,  $\phi^*$ is identical on $H^*(X^{\times n};\kk)$. This implies that
the differential $d_m$ of our spectral
sequence acts by $d_m(s_i) =q^\ast_i(\Delta)$, as stated in the theorem.
{\proofend}

\begin{remark}\label{rmk2} {\rm Without assumption $m>1$ the arguments of the
proof involving the top-dimensional cohomology classes do not
work. In fact, Theorem \ref{thm2}
fails for $X=S^1$ (see Example ~\ref{ex1}  above and  Remark
~\ref{rmk32}  below). The assumption $m>1$ should be also added to Theorem 1 of
Totaro \cite{To}.}
\end{remark}

\section{Cyclic configuration space of the sphere}

In this section we will use Theorem \ref{thm2}  to
calculate the cohomology ring of $G(S^m,n)$ with coefficients in $\Z_2$.

\begin{theorem}\label{thm3} Let $m>2$.
The cohomology ring $H^\ast(G(S^m,n);\Z_2)$ is multiplicatively
generated by cohomology classes
$$\sigma_i\in H^{i(m-1)}(G(S^m,n);\Z_2) \quad {\rm where} \quad i=1, 2,
\dots \quad
{\rm and}  \quad u\in H^{m}(G(S^m,n);\Z_2).$$
 These classes satisfy the following relations:
\begin{itemize}
\item[{\rm (i)}] $\sigma_i=0$ for $i\ge n-1$;
\item[{\rm (ii)}] $\sigma_i\sigma_j = \displaystyle{i+j \choose
i}\sigma_{i+j}$;
\item[{\rm (iii)}] $u^2=0$.
\end{itemize}
The Poincar\'e polynomial of the space $G(S^m,n)$ with coefficients in
$\Z_2$ equals
$${{(t^m+1) (t^{(n-1)(m-1)}-1)}\over{(t^{m-1}-1)}},$$
and the sum of Betti numbers is $2(n-1)$.
The dihedral group $D_n$ acts trivially on  $H^\ast(G(S^m,n);\Z_2)$.
\end{theorem}

\begin{remark}\label{rmk31} {\rm Relation (ii) implies $\sigma_i^2=0$, since
the binomial coefficient
$\displaystyle{2i\choose i}$ is always even.}\end{remark}

\begin{remark}\label{rmk32} {\rm  The space $G(S^1,n)$ consists of $n-1$
connected components and each is
homotopy equivalent
to $S^1$. Comparing with Theorem \ref{thm3},
we see that this theorem gives a correct additive structure of
$H^\ast(G(S^1,n);\Z_2)$,
although the multiplicative structure for $m=1$ is different. In fact,
for any zero-dimensional cohomology class $x$ with $\Z_2$ coefficients
one has $x^2=x$, and so for $m=1$ the classes $\sigma_i$ with
$\sigma_i^2=0$ do not exist. Nor is the action of the dihedral group $D_n$ on
$H^\ast(G(S^1,n);\Z_2)$ trivial
since the reflection changes the rotation number (as we mentioned above)
and so it acts nontrivially
on $H^0(G(S^1);\Z_2)$.}\end{remark}

\begin{remark}\label{rmk33} {\rm  Most of the statements of Theorem
\ref{thm3} hold true for $m=2$. In particular, for $m=2$, the
Poincar\'e polynomial
remains the same as stated in  Theorem \ref{thm3}. There
are only two
distinctions.
First, for $m=2$,  relation (ii) should be replaced by a more general one
\begin{itemize}
\item[(ii')]
$\sigma_i\sigma_j = \displaystyle{i+j \choose
i}\sigma_{i+j}+\varepsilon_{i,j}\sigma_{i+j-2}u$, where
$\varepsilon_{i,j}\in \Z_2$;
\end{itemize}  Secondly, our methods do
not prove that the
dihedral group $D_n$ acts trivially on
$H^\ast(G(S^2,n);\Z_2)$; instead we have a weaker statement that the cyclic
group $\Z_n\subset
D_n$ acts trivially on $H^\ast(G(S^2,n);\Z_2)$, if $n$ is odd.}
\end{remark}

\begin{remark}\label{rmk34} {\rm Paper \cite{Ba} contains an error.
Proposition 3.3 of \cite{Ba}
claims that $G(S^2,n)$ is simply connected and all further computations in
\cite{Ba}
depend on this claim.
In fact, the space $G(S^2,n)$ is simply connected only for $n=1$ and $n=2$,
as follows
from our previous remark \ref{rmk33}.
}\end{remark}

{\bf Proof.} In what follows we will omit the coefficients $\Z_2$ from the
notation.  We will apply the spectral sequence of Theorem \ref{thm2} in
the case when $X=S^m$. Let us describe its $E_2$-term.

Let $v\in H^m(S^m)$ be the generator. For $i=1, 2, \dots, n$
let $u_i\in H^m((S^m)^{\times n})$ denote the class
$u_i=1\times  \dots \times v\times 1\times \dots \times 1$ (with $v$ at the
$i$-th position). The classes $u_1, \dots, u_n$ generate the cohomology
ring $H^\ast((S^m)^{\times n})$; they commute and satisfy the relations
$u_i^2=0$.

The $E_2=E_m$-term of the spectral sequence  of Theorem \ref{thm2}
is the quotient of the bigraded commutative differential algebra
${\Z}_2 [s_1, \dots, s_n,u_1, \dots, u_n]$, where each
$s_i$ has bidegree $(0, m-1)$ and each $u_i$ has bidegree $(m,0)$,
by the ideal generated by the relations:
\begin{itemize}
\item[(a)] $s_i^2=0$ for $i=1, 2, \dots, n$;
\item[(b)]  $\sigma_{n-1} (s) = s_1s_2\dots s_{n-1}+s_2s_3\dots s_n +
s_3s_4\dots s_ns_1+\dots
+s_ns_1\dots s_{n-2}=0$;
\item[(c)]  $u_i^2=0$ for $i=1, 2, \dots, n$;
\item[(d)]  $(u_i+u_{i+1})s_i =0$ for $i=1, 2, \dots, n$.
\end{itemize}

Recall that we understand indices cyclically; for example, (d) contains the
relation $(u_1+u_n)s_n=0$.

The first nontrivial differential is $d_m$ and it acts as follows
\[d_m s_i = u_i + u_{i+1},\ d_m u_i =0, \quad i=1,...,n.\]

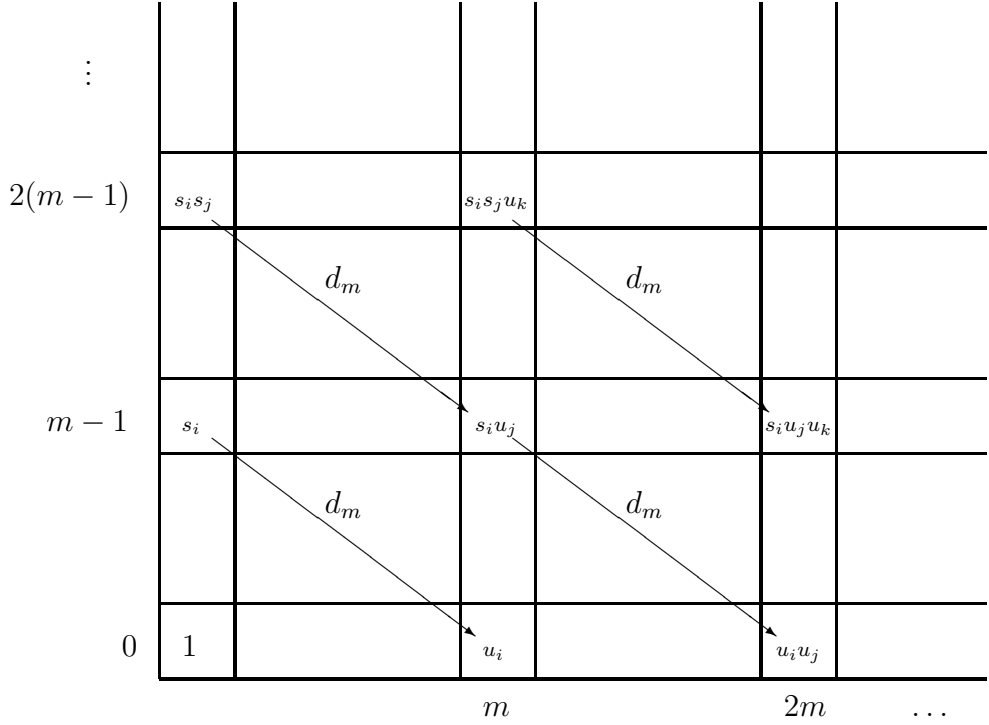
\begin{figure}[h]
\setlength{\unitlength}{1cm}
\begin{center}
\begin{picture}(12,10)
\linethickness{0.3mm}
\multiput(0,0)(0,3){3}{\line(1,0){11}}
\multiput(0,1)(0,3){3}{\line(1,0){11}}
\multiput(0,0)(4,0){3}{\line(0,1){9}}
\multiput(1,0)(4,0){3}{\line(0,1){9}}
\put(0.3,0.3){1}
\put(4.3, 0.3){$\scriptstyle u_i$}
\put(8.2, 0.3){$\scriptstyle u_iu_j$}
\put(8.05,3.3){$\scriptstyle s_iu_ju_k$}
\put(0.3, 3.3){$\scriptstyle s_i$}
\put(4.2, 3.3){$\scriptstyle s_iu_j$}
\put(0.2, 6.3){$\scriptstyle s_is_j$}
\put(4.05, 6.3){$\scriptstyle s_is_ju_k$}
\put(4.3,-0.5){$m$}
\put(8.3, -0.5){$2m$}
\put(10, -0.5){$\dots$}
\put(-0.5,0.3){$0$}
\put(-1.5, 3.3){$m-1$}
\put(-2, 6.3){$2(m-1)$}
\put(-1, 7.9){$\vdots$}
\put(0.7,3.2){\vector(4, -3){3.5}}
\put(0.7,6.1){\vector(4, -3){3.4}}
\put(2.2,2.2){$d_m$}
\put(2.2,5.2){$d_m$}
\put(4.7,3.2){\vector(4, -3){3.5}}
\put(4.7,6.1){\vector(4, -3){3.4}}
\put(6.2,2.2){$d_m$}
\put(6.2,5.2){$d_m$}
\end{picture}
\end{center}
\caption{Term $E_m$ of the spectral sequence}
\end{figure}

Introduce new variables: $v_i =  u_i + u_{i+1}$, where $i=1, 2, \dots, n$
and $u = u_1$.
In the new variables the differential algebra $(E_m,d_m)$ can be understood
as the quotient of
the polynomial ring
$\Z_2[s_1, \dots, s_n, v_1, \dots, v_n, u]$ modulo the ideal generated by
the relations
(a), (b) and the following relations (c'), (d') and (e):
\begin{itemize}
\item[(c')]  $v_i^2=0$ for $i=1, 2, \dots, n$ and $u^2=0$;
\item[(d')]  $v_is_i =0$ for $i=1, 2, \dots, n$;
\item[(e)]  $\sum_{i=1}^n v_i=0$.
\end{itemize}
The differential $d_m$ acts by
\begin{eqnarray}
d_m(s_i) = v_i, \quad  d_m(v_i)=d_m(u)=0.\label{diff}
\end{eqnarray}
Our goal is to compute the cohomology of $(E_m, d_m)$.

Let $A$ denote the quotient of $\Z_2[s_1, \dots, s_n, v_1, \dots, v_n, u]$
by the ideal
generated by the relations (a), (c'), (d'), i.e., we simply ignore
relations (b) and (e).
We consider $A$ and $E_m$
with the total grading, where each $s_i$ has degree $m-1$ and the elements
$v_i$ and $u$ have degree $m$. Let the differential
$d_A: A\to A$ act by $d_A(s_i)=v_i$ and $d_A(v_i)=0=d_A(u)$. There is a
canonical
epimorphism of graded differential algebras $f:A\to E_m$, so that $fd_A=d_m
f$.
The kernel of $f$ is the ideal of $A$ generated by $\sum_{i=1}^n v_i$ and
$\sigma_{n-1}(s)$.

We claim that
\begin{eqnarray}
H^i(A,d_A)=\left\{
\begin{array}{l}
0, \quad {\rm if}\quad i\ne m,\\ \\
\Z_2,\quad {\rm if}\quad i=m,
\end{array}\label{acyc1}
\right.
\end{eqnarray}
and the only nontrivial cohomology class is represented by $u$. In order to
prove this, consider the
filtration $A_0\subset A_1\subset \dots \subset A_n=A$, where $A_0$ is
generated by $u$,
and each $A_k$ is generated by $u, s_i, v_i$ for $i=1, \dots, k$.
The differential $d_A$ preserves this filtration.
$d_A$, restricted to $A_0$, vanishes, and so
$(A_0,d_A)$ has a one-dimensional cohomology generated by the class of $u$.
We will now show that all factors $A_k/A_{k-1}$ are acyclic, where $k=1,
\dots, n$.
Indeed, any element $a\in A_k/A_{k-1}$ can be uniquely represented in the form
$a=s_kx +v_k y$, where $x, y\in A_{k-1}$. If $a$ is a cocycle then we have
$$d_A(a)= v_k x+s_kd_A(x) +v_kd_A(y)=s_kd_A(x)+v_k[x+d_A(y)]=0$$
 and hence $x=d_A(y)$. Therefore $a=d_A(s_ky)$. This proves that
$H^\ast(A,d_A)=H^\ast(A_0,d_A)$
and our claim (\ref{acyc1}) follows.

Introduce a new differential $\delta_A: A\to A$ of degree $m$:
$$\delta_A(x)=(\sum_{i=1}^n v_i)x.$$
Clearly, $\delta_A^2=0$ and $\delta_A d_A =d_A\delta_A$; however $\delta_A$
does not obey the
Leibnitz rule. We claim that
\begin{eqnarray}
H^i(A,\delta_A)=\left\{
\begin{array}{l}
0, \quad {\rm if}\quad i\ne n(m-1), \quad i\ne n(m-1)+m,\\ \\
\Z_2,\quad {\rm if}\quad i=n(m-1)\quad {\rm or}\quad i=n(m-1)+m,
\end{array}
\right.\label{acyc2}
\end{eqnarray}
and the only nontrivial cohomology classes are represented by $s_1s_2\dots
s_n$ and by
$s_1s_2\dots s_nu$.
Indeed, each element of $A$ can be written as a sum of monomials in $s_i$,
$v_i$ and $u$.
For $I\subset \{1, 2, \dots n\}$,  denote by $s_I$ the product of all $s_i$
for $i\in I$. Similarly,  label monomials in variables $v_i$ as $v_J$ for
$J\subset \{1, 2, \dots n\}$. Note that the product
$s_I v_J\in A$ is nontrivial if and only if $I$ and $J$ are disjoint
subsets of $\{1, \dots, n\}$.
Note also that
\[\delta_A(s_Iv_J)\, =\,  \sum_{i\notin I\cup J} s_Iv_{J\cup \{i\}}, \quad
\delta_A(s_Iv_Ju)\, =\,  \sum_{i\notin I\cup J} s_Iv_{J\cup \{i\}}u.\]
We see that application of $\delta_A$ does not change the multi-index $I$
or the factor $u$.
Hence, the complex $(A, \delta_A)$ splits into a direct sum over different
multi-indices $I$.
Fix a set $I$ and denote by $k$ the cardinality of the set $\{1,\dots,n\} -
I$.
Then the respective part of the complex $(A, \delta_A)$ is isomorphic to
two copies
(one with $u$ and one without)
of the standard cochain complex of the simplex with $k$ vertices:
the differential of an  $r$-dimensional face (i.e., set $J$) is the sum of
$r+1$-dimensional faces
that contain the given one (sets $J \cup \{i\}$). Note that empty set $J$
is also allowed.
This complex has zero cohomology, unless $k=0$
(empty simplex), in which case the cohomology is ${\Z}_2$.
This exceptional case corresponds to $I = \{1,\dots,n\}$, and the claim
follows.

Now we are ready to compute the cohomology of $(E_m, d_m)$,
i.e., the term $E_{m+1}$ of the spectral sequence.
Denote by $\sigma_k\in A$ the $k$-th elementary symmetric function in
variables
$s_1\dots, s_n$, i.e.,
\[\sigma_0=1, \quad{\rm and}\quad
\sigma_k =\sum_{1\le i_1<\dots <i_k\le n} s_{i_1}s_{i_2}\dots
s_{i_k}\quad{\rm for}\quad
 k=1, 2, \dots, n.\]
It is clear that
$$
\sigma_i\sigma_j \, =\, \left\{
\begin{array}{l}
\displaystyle {i+j \choose i}\sigma_{i+1}, \quad{\rm for}\quad i+j\le n,\\ \\
0,\quad{\rm for}\quad i+j>n.
\end{array}
\right.
$$
Also note that
\begin{eqnarray}
d_A(\sigma_{k+1}) = \delta_A(\sigma_k), \quad
d_A(\sigma_{k+1}u) = \delta_A(\sigma_k u)\label{sigma}.
\end{eqnarray}
Therefore the images of the classes $\sigma_i$ and $\sigma_i u$ under the
projection
$f: A\to E_m$ are cocycles.
Note also that the classes $f(\sigma_{n})$,
$f(\sigma_{n-1})$ vanish due to the relation (b).
All other classes $f(\sigma_0)=1, f(\sigma_1), \dots, f(\sigma_{n-2})$
are nonzero elements of $E_m$;
this follows since the
kernel of $f$ is generated by the image of $\delta_A$ and by the element
$\sigma_{n-1}$, cf. above.
 It is clear from dimension considerations that $\sigma_i$ is not a
coboundary; likewise,
$\sigma_i u$ is not a coboundary either.
Hence we have  nontrivial cohomology
classes
$$f(\sigma_i)\in H^{i(m-1)}(E_m,d_m) \quad {\rm and} \quad f(\sigma_i u)\in
H^{i(m-1)+m}(E_m,d_m), \quad  i=0, 1, \dots, n-2.$$
Our aim now is to show that these classes constitute all the cohomology.

Let $K\subset A$ denote the kernel of $f$.  Using long exact sequences, we
deduce from (\ref{acyc1})
 that
\begin{eqnarray}
H^i(E_m,d_m) \, \simeq\,  H^{i+1}(K,d_A)\quad{\rm  for}\quad i\ne m.
\label{shift}
\end{eqnarray}
Note that $K$ contains $\sigma_{n}$, and  denote by $\overline K$ the
factor of $K$
by the ideal, generated by $\sigma_{n}$ and $d_A(\sigma_n)$. This ideal is
four-dimensional and is
 generated, as a vector space, by the classes $\sigma_n$, $d_A(\sigma_n)$,
$\sigma_nu$ and
$d_A(\sigma_n)u$; it is closed under $d_A$ and is acyclic.  Hence we have
$$H^\ast(K,d_A) \, \simeq\,  H^\ast(\overline K, d_A).$$
Define a
chain map $\phi: (E_m,d_m) \to (\overline K, d_A)$ increasing the degree by
$m$. Given an
element $x\in E_m$, consider a lift $\overline x\in A$  and set:
$\phi(x)\, =\, \delta_A(\overline x) \, \in \, \overline K.$
To see that $\phi$ is well-defined, note that $\overline x$ is determined
up to summation with
elements of the form $\delta_A(a)+\sigma_{n-1}b$, where $a, b\in A$. Since
$\delta_A(\delta_A(a) +\sigma_{n-1} b)\, =\, d_A(\sigma_n) b \, =\, 0 \in
\overline K$,
we see that $\phi(x)$ does not depend on the choice of $\overline x$.
It is clear that $\phi$ is a chain map, i.e.,
$\phi d_m = d_A \phi$. Due to (\ref{acyc2}), $\phi$ is a monomorphism.
The cokernel of $\phi$ is linearly generated by the classes $\sigma_{n-1}$
and $\sigma_{n-1}u$,
and we find from the short exact sequence
$0\to E_m \stackrel {\phi}{\to} \overline K \to \coker(\phi)\to 0$
 that
\begin{eqnarray}
\phi: H^i(E_m, d_m) \, \to\,  H^{i+m}(\overline K, d_A)\label{period}
\end{eqnarray}
is an isomorphism for all $i$ except $i=(n-2)(m-1)$ and $i=(n-1)(m-1) +1$;
for these two exceptional
values of $i$ the homomorphism (\ref{period}) is an epimorphism with one
dimensional kernel.
To prove this claim about (\ref{period}) notice that the induced
differential on
$\coker(\phi)$ is trivial, as well as the induced map
$H^\ast(\overline K) \to H^\ast(\coker(\phi))$.

Combining isomorphisms (\ref{shift}) and (\ref{period}), one finds that
\begin{eqnarray}
H^i(E_m,d_m) \, \, \stackrel{\simeq}{\to} \, \,  H^{i+m-1}(E_m,
d_m),\label{period1}
\end{eqnarray}
where $i\ne 1$, $i\ne (n-2)(m-1)$ and $i\ne (n-1)(m-1) +1$. Note also that for
$i = (n-2)(m-1)$ and $i =(n-1)(m-1) +1$ instead of isomorphism
(\ref{period1}) we have an
epimorphism with one dimensional kernel.

\begin{figure}[h]
\setlength{\unitlength}{0.9cm}
\begin{center}
\begin{picture}(10,15)
\linethickness{0.35mm}
\multiput(0,0)(0,3){3}{\line(1,0){9}}
\multiput(0,1)(0,3){3}{\line(1,0){9}}
\put(0,12){\line(1,0){9}}
\put(0,13){\line(1,0){9}}
\multiput(0,0)(4,0){2}{\line(0,1){14}}
\multiput(1,0)(4,0){2}{\line(0,1){14}}
\put(0.3, 0.3){$\scriptstyle 1$}
\put(4.3, -0.5){$m$}
\put(4.3, 0.3){$\scriptstyle u$}
\put(4.2, 3.3){$\scriptstyle \sigma_1 u$}
\put(4.2, 6.3){$\scriptstyle \sigma_2 u$}
\put(0.3, 3.3){$\scriptstyle \sigma_1$}
\put(0.3, 6.3){$\scriptstyle \sigma_2$}
\put(0.1, 12.3){$\scriptstyle \sigma_{n-2}$}
\put(4.05, 12.3){$\scriptstyle \sigma_{n-2}u$}
\multiput(0.5, 8)(2,0){5}{$\vdots$}
\multiput(0.5, 10)(2,0){5}{$\vdots$}
\put(-1.5, 3.3){$m-1$}
\put(-1.8, 6.3){$2m-2$}
\put(-3.5,12.3){$(n-2)(m-1)$}
\end{picture}
\end{center}
\caption{Term $E_{m+1}$ of the spectral sequence.}
\end{figure}

 We can finally show that the term $E_{m+1}$ is linearly
generated by the classes $\sigma_1, \dots, \sigma_{n-2}$ and $u, \sigma_1u,
\dots,
\sigma_{n-2}u$, where, abusing notation, we identify $f(\sigma_i)$ and
$f(u)$ with $\sigma_i$
and $u$, respectively.
Obviously, $\sigma_1$ and $u$ are the only nontrivial classes of
$E_{m+1}$ of total degree $<2(m-1)$. The periodicity isomorphism
(\ref{period1}) implies that in
each degree $i(m-1)$ (where $1\le i<n-1$) we have only one nontrivial
class, which therefore must
coincide with
$\sigma_i$. Also, again from (\ref{period1}), one concludes that in each
degree $i(m-1)+m$, where
$1\le i<n-1$, there is a single nontrivial class, and so it must coincide
with $\sigma_iu$.
Hence the term  $E_{m+1}$ has the structure shown in Figure 2.

By dimension considerations, all further differentials
$d_r$ with $r>m$ vanish and thus $E_{m+1}=E_{\infty}$.

Now we want to reconstruct the cohomology algebra $H^\ast(G(S^m,n))$ from
$E_\infty$.
The $m$-th column  is an ideal in $H^\ast(G(S^m,n))$, and the factor of
$H^\ast(G(S^m,n))$ by this
ideal is the $0$-th column. Since $m>2$,
each diagonal $x+y=c$ contains at most one nonzero term of $E_\infty$.
Hence the generators $\sigma_i$ and $u$ admit unique lifts to
$H^*(G(S^m,n))$, and  we  label
the cohomology classes of $H^\ast(G(S^m,n))$ as $\sigma_i$ and $\sigma_iu$.
Since the multiplication in $E_\infty$ is induced from the multiplication
in $H^\ast(G(S^m,n))$,
 the product of $\sigma_iu$ and $\sigma_ju$ is trivial, and the product of
$\sigma_i$ and
$\sigma_ju$ equals $ {i+j\choose i}\sigma_{i+j}u$.
The product of $\sigma_i$ and $\sigma_j$ may
be equal to ${i+j\choose i}\sigma_{i+j}$ plus a term from the $m$-th column,
but, in fact,  this additional term
 must vanish  since the $m$-th column has no nonzero terms
in  dimensions, divisible by $m-1$. We conclude
that the cohomology
algebra $H^\ast(G(S^m,n))$ has the structure described in Theorem \ref{thm3}.

The dihedral group $D_n$ acts trivially on $H^\ast(G(S^m,n))$
since in each dimension we have at most one nonzero element.

The arguments in the two last paragraphs do not apply if $m=2$.
This explains the modifications to the statement of Theorem \ref{thm3} in
case $m=2$, described
in Remark ~\ref{rmk33} above. Let us show that $Z_n\subset D_n$
acts trivially on $H^\ast(G(S^2,n);\Z_2)$, assuming that $n$ is odd.
Indeed, it follows from the above calculation of the $E_\infty$-term that
$H^i(G(S^2,n);\Z_2)$ is at
most two-dimensional for any $i$, and hence it has either 1 or 3 nonzero
elements.
$\Z_n$ acts by permutations of nonzero classes, and at least
one of these classes is fixed, since $D_n$ acts trivially on the
$E_\infty$-term.
Therefore the action of $\Z_n$ is trivial provided $n$ is odd. This
completes the proof.
\proofend

\begin{corollary}\label{cor1}
Consider the cohomology classes
\[\sigma_1, \, \sigma_2, \, \sigma_4, \dots, \sigma_{2^k} \, \in\,
H^\ast(G(S^m,n);\Z_2),
 \quad\quad m>2,\]
where $k$ is the largest integer,  not exceeding $\log_2(n-2)$.
Then, for any $i \leq n-2$, the class $\sigma_i$ can be expressed as a product
\begin{eqnarray}
\sigma_i = \sigma_{2^{k_1}} \sigma_{2^{k_2}} \dots \sigma_{2^{k_r}},
\end{eqnarray}
where $i = 2^{k_1} + 2^{k_2}+ \dots +2^{k_r}$ is the binary expansion of $i$;
here $k_1 < k_2< \dots <k_r$.
\end{corollary}

{\bf Proof.} It follows from  relation (ii) of Theorem \ref{thm3} that
$$\sigma_{l_1} \dots \sigma_{l_r} = {l_1 + \dots + l_r \choose l_1, \dots,
l_r}
\sigma_{l_1 + \dots +l_r},\quad \quad
l_1 + \dots + l_r  \le n-2.$$
Hence the statement of the Corollary will follow once we show that the
multinomial coefficient
$${2^{k_1} + \dots + 2^{k_r} \choose 2^{k_1}, \dots, 2^{k_r}}=
{{(2^{k_1} + \dots + 2^{k_r})!}\over({2^{k_1})! \dots (2^{k_r})!}}$$
is odd for $k_1<\dots <k_r$. This is a consequence of the well known fact that
the maximal power of $2$ that divides $l!$ equals $l-r$, where $r$ is the
number of
1's in the binary expansion of $l$. {\proofend}

\begin{corollary}\label{cor2}
For $m>1$ the  cup-length of $H^\ast(G(S^m,n);\Z_2)$ equals
\[ [\log_2(n-1)]+1.\]
\end{corollary}
{\bf Proof.} Assume first that $m>2$.
Then, by Corollary \ref{cor1}, the length of the longest nontrivial product
of classes $\sigma_i$
equals the maximal number of 1's in the binary expansion, that  integers,
not exceeding $n-2$, may have. It is easy to see that the latter number
equals $[\log_2(n-1)]$.
This implies Corollary \ref{cor2},
since we  also have the class $u$ at our disposal.
For $m=2$ the arguments are similar using relation (ii') in Remark
~\ref{rmk33}.
\proofend

\section{Morse theory of closed billiard trajectories}

In this section we describe  Morse theory of closed billiard
trajectories
which is similar to the classical Morse theory of closed geodesics in
Riemannian manifolds. Some of the results of this section are known to experts,
but do not exist in the literature in an accessible form.

The proof of Theorem \ref{thm1} will depend on the results of this section
in two ways. First,  Morse theory of closed billiard trajectories, in the
simplest case
when the billiard table
is the round sphere (like in \cite{Ba}), will provide a tool for
computation of the
$D_n$-equivariant cohomology
ring of the cyclic configuration space $G(S^m,n)$. Second,
we will invoke  equivariant Morse theory \cite{Bo}
and a version of equivariant Lusternik-Schnirelman theory developed in
\cite{Mar}  to  deduce Theorem \ref{thm1} concerning arbitrary smooth strictly
convex billiards.

Let $X \subset {\R}^{m+1}$ be a smooth closed strictly
convex hypersurface, topologically the sphere, which is the boundary
of the billiard table. Denote by
\begin{eqnarray}
L_X: G(X,n)\to \R
\end{eqnarray}\label{length}
the perimeter length function, taken with the minus sign,
\begin{eqnarray*}
L_X(x_1,\dots,x_n) = -|x_1 - x_2| - |x_2 - x_3| - \dots - |x_n - x_1|,
\end{eqnarray*}
where $(x_1, x_2, \dots, x_n)\in G(X,n)$ and the distance $|x_i-x_{i+1}|$
is measured in the ambient Euclidean space $\R^{m+1}$. The reason for the
minus sign
will become clear shortly.

It is well-known that
$n$-periodic billiard orbits in $X$ are precisely the critical points of
the function $L_X$; this is the Maupertuis'
principle of the classical mechanics as applied to billiards (see, e.g.,
\cite {KT,Ta}).
Clearly, the function $L_X$ is smooth and
$D_n$-equivariant. Identifying $G(X,n)$ with $G(S^m, n)$, we see that the
shape of the billiard
domain $X$ becomes encoded in the function $L_X: G(S^m,n)\to \R$, and the
problem of finding
the closed billiard trajectories inside $X$ turns into a problem of
Morse-Lusternik-Schnirelman theory.

We encounter the following difficulty: one cannot apply
Morse-Lusternik-Schni\-relman
theory directly to $G(X,n)$ since
this manifold is not compact. The function $L_X$  extends to a continuous
function on
the space of all $n$-tuples $X^{\times n}$
but this extension fails to be differentiable on the singular set
$\Sigma$ consisting of the
points with $x_i = x_{i+1}$ for some $i$. A way around this difficulty is
in replacing
$G(X,n)$ by a compact manifold with boundary $G_{\varepsilon} (X,n) \subset
G(X,n)$, where $\varepsilon >0$ is small enough and
\begin{eqnarray}
G_{\varepsilon} (X,n) = \{(x_1, \dots, x_n) \in X^{\times n} :
\prod_{i=1}^n|x_i - x_{i+1}| \geq \varepsilon \};
\end{eqnarray}
similar approach can be found in \cite {Ba} and in \cite {CS}, \cite {KT}
for the two-dimensional
case.

\begin{proposition} \label{prop2} If $\varepsilon>0$ is sufficiently small
then:
\begin{enumerate}
\item[(a)] $G_{\varepsilon} (X,n)$ is a smooth manifold with boundary;
\item[(b)] the inclusion $G_{\varepsilon} (X,n)\subset G(X, n)$ is a
$D_n$-equivariant homotopy equivalence;
\item[(c)] all
critical points of $L_X: G(X, n)\to \R$ are contained in $G_{\varepsilon}
(X,n)$;
\item[(d)] at every point of  $\partial G_{\varepsilon} (X,n)$,
the gradient of $L_X$ has the outward direction.
\end{enumerate}
\end{proposition}

{\bf Proof}. The function $\phi (x_1, \dots, x_n) = \prod_1^n |x_i -
x_{i+1}|^2$ is
smooth on $X^{\times n}$. Its zero level set $\phi^{-1}(0)$
is $\Sigma = X^{\times n} - G(X,n)$, which is a
critical level. Statements (a) and (b) will follow once we show that {\it there
exists a constant $\delta>0$ such that the interval $(0, \delta)$
consists of regular values of $\phi$.}

To prove the claim in italic we need to understand geometry of $n$-tuples
$\bar x = (x_1, \dots, x_n)\in G(X,n)$
which are critical points of $\phi$. An easy computation with Lagrange
multipliers shows that $\bar x$ is critical if and only if
\begin{eqnarray}
\frac{x_{i-1} - x_{i}}{|x_{i-1}- x_{i}|^2} +
\frac{x_{i+1}- x_{i}}{|x_{i+1}-x_{i}|^2} = t_i \nu_i,
\quad i=1,\dots,n \label{crit}
\end{eqnarray}
where $\nu_i = \nu (x_i)$ is the unit normal vector to $X$ at point $x_i$
and $t_i\in \R$
is some constant. It will be convenient for us to assume that $\nu_i$ has
the inward direction.

Since $X$ is smooth and strictly convex one can find two positive numbers
$r < R$ such
that for every  $x \in X$

(1) there is a sphere $s(x)$ of radius $r$, tangent to $X$ at $x$ and
contained inside $X$; and

(2) there is a sphere $S(x)$ of radius $R$,
tangent to $X$ at $x$ and containing $X$.

Let $x_{i-1}, x_i, x_{i+1}$
be three consecutive points on $X$ from a critical $n$-tuple $\bar x$ of the
function $\phi$.  Note that the
2-plane through the points $x_{i-1}, x_i, x_{i+1}$ contains the centers of
both spheres $s(x_i)$ and $S(x_i)$.  Let
$\alpha_i$ be the angle between $\nu_i$ and $x_{i-1}-x_i$, and $\beta_i$
the angle between $\nu_i$ and $x_{i+1}-x_i$. Set: $l_i=|x_i - x_{i-1}|$.
Then (\ref{crit}) gives
\begin{eqnarray}
{\sin \alpha_i}/{l_i} ={\sin \beta_i}/{l_{i+1}}.\label{crit2}
\end{eqnarray}
Since $x_{i-1}$ and $x_{i+1}$ lie outside the sphere $s(x_i)$ and inside
the sphere $S(x_i)$,
we have
\begin{eqnarray}
2r \cos \alpha_i < l_i < 2R \cos \alpha_i,\quad 2r \cos \beta_i < l_{i+1} <
2R \cos \beta_i, \label{ineq}
\end{eqnarray}
and hence
\begin{eqnarray}
 {{r\cos \alpha_i}\over{R\cos \beta_i}} \, \leq {{l_i}\over{l_{i+1}}}\, \leq
{{R\cos \alpha_i}\over{r\cos \beta_i}}.\label{pin1}
\end{eqnarray}
We claim that
\begin{eqnarray}
{{r}\over{\sqrt{r^2 + R^2}}} \, \leq {{l_i}\over{l_{i+1}}} \leq \,
{{\sqrt{r^2 +
R^2}}\over{r}} \label{pin}
\end{eqnarray}
for any $i=1, 2, \dots, n$. Indeed,
if (\ref{pin}) fails and
\begin{eqnarray}
{{l_i}\over{l_{i+1}}} < {{r}\over{\sqrt{r^2 + R^2}}} \label{if}
\end{eqnarray}
then, combining (\ref{if}) with the left inequality in (\ref{pin1}), we obtain
\begin{eqnarray}
\cos \alpha_i < {{R\cos \beta_i}\over{\sqrt{r^2 + R^2}}} \leq
{{R}\over{\sqrt{r^2 +
R^2}}}. \label{cos}
\end{eqnarray}
On the other hand, combining (\ref{if}) with (\ref{crit2}), we obtain
\begin{eqnarray}
\sin \alpha_i < {{r \sin \beta_i}\over{\sqrt{r^2 + R^2}}} \leq
{{r}\over{\sqrt{r^2 + R^2}}}, \label{sin}
\end{eqnarray}
which leads to a contradiction since (\ref{cos}) and (\ref{sin}) are
incompatible.
This argument shows that the left inequality in (\ref{pin}) holds.
The right inequality in (\ref{pin}) follows similarly.

It follows from (\ref{pin})  that for any two edges $l_i$ and $l_j$ of a
critical
$n$-gon one has:
\begin{eqnarray*}
{{l_i}\over{l_j}} \, <\, \Bigl(1+{{R^2}\over{r^2}}\Bigr)^{n/2}\label{pin2}.
\end{eqnarray*}
In other words, if one of the edges of a critical $n$-gon is "short" then
all others are
also "short".

It remains to show that a critical $n$-gon cannot have all edges
arbitrarily short.
Indeed, assume that $\bar x= (x_1,\dots,x_n)$ is a critical polygon such that
all points $x_i$ lie in a small neighborhood  $U \subset X$ of a point
$y\in U$.
We will assume that
$U$ is so small that for any $x\in U$ the scalar product $\langle \nu(x),
\nu(y)\rangle$ is
positive
(recall that $\nu(x)$ denotes the inward unit normal to $X$ at point $x$).
Taking scalar product of equation (\ref{crit}) with $\nu_i$ gives
$$
{{\cos \alpha_i}\over{l_i}} + {{\cos \beta_i}\over{l_{i+1}}} =t_i,
$$
and hence  all numbers $t_i$, which appear in (\ref{crit}), are positive.
Now, the scalar product of (\ref{crit}) with vector $\nu(y)$ gives
\begin{eqnarray}
{{\langle ({x_{i-1}-x_i}), \nu(y)\rangle}\over{l_i^2}} \, >\, {{\langle
({x_{i}-x_{i+1}}),
\nu(y)\rangle}\over{{l_{i+1}^2}}},\quad i=1, 2, \dots, n.\label{cyclic1}
\end{eqnarray}
Recall that in the last inequality the indices are understood cyclically,
i.e. $x_{n+j}=x_j$.  Hence  (\ref{cyclic1}) leads to a contradiction.

The above arguments prove statements (a) and (b).

Next we prove statement (c). The argument is similar. If $\bar x = (x_1,
\dots, x_n)\in G(X,n)$ is
a billiard trajectory in $X$, then (instead of (\ref{crit})) we have
\begin{eqnarray*}
{{x_{i-1} - x_{i}}\over{l_i}} +
{{x_{i+1}- x_{i}}\over{l_{i+1}}} = t_i \nu_i,
\quad i=1,\dots,n, \label{crit3}
\end{eqnarray*}
and (\ref{crit2}) becomes the usual reflection law
$$\alpha_i =\beta_i.$$
We have the inequality
$${{r}\over{R}} \, \leq \, {{l_i}\over{l_{i+1}}} \leq {{R}\over{r}}$$
which is an analog of (\ref{pin}) in the present case.
As above, this inequality implies that if one edge of an $n$-periodic
billiard trajectory is "short"
then so are all its edges. The preceding argument  shows that billiard $n$-gons
cannot lie entirely in a small neighborhood of a point of the hypersurface $X$.

Finally we prove claim (d).
Since the gradient $\nabla \phi$ is orthogonal to the
boundary $\partial G_\varepsilon(X,n)$ and points inside $G_\varepsilon(X,n)$,
it suffices to show that the scalar product
$\langle \nabla (\ln \phi), \nabla L_X \rangle$ is negative along
$\partial G_\varepsilon(X,n)$ for every sufficiently small $\varepsilon >0$.

Let us compute the gradients involved.
Taking into account the decomposition
$T_{\bar x} G(X,n) = T_{x_1} X \times \dots \times T_{x_n} X,$ where $\bar
x\in G(X,n)$,
one finds that the $i$-th components of the
gradients $\nabla L_X$ and $\nabla (\ln \phi)$ are given by
\begin{eqnarray}
(\nabla L_X)_i = -{{{x_{i-1}-x_i}}\over{l_i}} -
{{{x_{i+1}-x_i}}\over{l_{i+1}}} + (\cos \alpha_i + \cos
\beta_i) \nu_i, \label{grad2}
\end{eqnarray}
\begin{eqnarray}
{{1}\over{2}}\ (\nabla \ln \phi)_i = {{{{{x_{i-1}-x_i}}\over{l_i^2}}}} +
{{{x_{i+1}-x_i}}\over{l_{i+1}^2}} - \Bigl({{\cos
\alpha_i}\over{l_i}} + {{\cos \beta_i}\over{l_{i+1}}}\Bigr) \nu_i.
\label{grad3}
\end{eqnarray}
Denote by $\theta_i$ the angle
between the vectors $x_{i-1}-x_i$ and $x_{i+1}-x_i$. Due to strict
convexity of $X$,
we have $\theta_i\in [0,\pi)$.
A direct computation using (\ref{grad2}) and (\ref{grad3}) shows that
\begin{eqnarray}
 {\langle \nabla L, \nabla \ln \phi \rangle} =  2( S_1 +S_2),\label{dot}
\end{eqnarray}
where
\begin{eqnarray}
S_1 = -\sum_{i=1}^n
\Bigl({{1}\over{l_i}} + {{1}\over{l_{i+1}}}\Bigr) (1 + \cos
\theta_i), \quad  S_2 = \sum_{i=1}^n (\cos \alpha_i + \cos \beta_i)
\Bigl({{\cos
\alpha_i}\over{l_i}} + {{\cos \beta_i}\over{l_{i+1}}}\Bigr). \label{dot1}
\end{eqnarray}

We want to show that the right-hand side of (\ref{dot}) is negative for
all $\bar x\in \partial G_\varepsilon (X,n)$
with $\varepsilon>0$
small enough.
It follows from inequalities (\ref{ineq}) that
$${{\cos \alpha_i}\over{l_i}} + {{\cos \beta_i}\over{l_{i+1}}} <
{{1}\over{r}},$$
therefore $S_2 < 2n/r$. We will be done once we show that $S_1$ tends to
$-\infty$
as $\varepsilon$ goes to zero.

Assume, to the contrary, that there exists a constant $C$ and an infinite
sequence
$\bar x_k\in G(X,n)$, where $k=1, 2, \dots$, such that
$\phi(\bar x_k) = \varepsilon_k$ tends to $0$ and
$S_1(\bar x_k) >-C$.
One has
\begin{eqnarray}
 \Bigl({{1}\over{l_i}} + {{1}\over{l_{i+1}}}\Bigr) (1 + \cos
\theta_i) < C, \quad i=1,\dots,n, \label{inq}
\end{eqnarray}
and $l_1 \dots l_n = \varepsilon_k^{1/2}$.

Suppose that for an $n$-gon $\bar x_k\in \partial G_{\varepsilon_k}(X,n)$
and some index $j$ one has an inequality
$l_i \le b\epsilon_k^a,$
where $k$ is large enough, $a>0$ and $b>0$ (for example, for the smallest
link $l_i$, one has: $l_i\le (\varepsilon_k)^{1/2n}$.
Then  it follows from (\ref{inq}) that
$$1 + \cos \theta_j < Cb\varepsilon_k^a,\quad {\rm and \, \, hence}\quad
\pi - \theta_j < {2\over \sqrt{3}} \sqrt{Cb}\varepsilon_k^{a/2}.$$
Since $\theta_j \leq \alpha_j + \beta_j$,  one concludes that
$$\pi/2 - \beta_j < {2\over \sqrt{3}} \sqrt{Cb}\varepsilon_k^{a/2},\quad
{\rm and}\quad
\cos \beta_j = \sin({\pi\over 2}-\beta_i) < {2\over \sqrt{3}}
\sqrt{Cb}\varepsilon_k^{a/2}.$$
It follows then from inequalities (\ref{ineq}) that
\begin{eqnarray}
l_{j+1} \, <\,  {4R\over \sqrt{3}} \sqrt{Cb}\varepsilon_k^{a/2}. \label{est2}
\end{eqnarray}

The argument that derives (\ref{est2}) from the initial assumption $l_i \le
b\epsilon_k^a$
can be repeated $n$ times
to conclude that there exist positive constants $b, b'$ and $a_0$,
so that for any $j=1, 2, \dots, n$ and for any
sufficiently large $k$ one has:
$$l_j \, <\,  b\varepsilon_k^{a_0} \quad {\rm and} \quad \pi - \theta_j \, <\,
b'\varepsilon_k^{{{a_0}/2}}$$
(we may take $a_0= 2^{-n}n^{-1}$).

Thus, for $k$ large enough, $\sum_{i=1}^n (\pi - \theta_i)$ is close to
zero. On the other
hand, given a closed polygonal line in Euclidean space, the sum of its
exterior angles, that
is,  the angles $\pi - \theta_i$, is at least $2\pi$ (a smooth version of
this statement holds
too: the total curvature of a closed curve is at least $2\pi$). This is a
contradiction.
\proofend

\begin{remark}\label{rmk4} {\rm $L_X$ and $\phi$ are particular
cases of the following more general
function on inscribed polygons:
$$F(x_1,\dots,x_n) = \sum_1^n f(|x_{i+1} - x_i|)$$
where $f(t)$ is a function of one variable (we obtain $L_X$ when $f(t) = t$
and $\ln \phi$ when $f(t) =
2\ln t$). For some of these functions it is not true  that all the critical
$n$-tuples lie off a neighborhood of the singular set $\Sigma$.  The
simplest example is provided by $f(t) = t^2$. An
analog of the condition (\ref{crit}) reads
$$\overrightarrow {x_i x_{i-1}} + \overrightarrow {x_i x_{i+1}} = t_i \nu
(x_i).$$
If $X$ is a circle then this equation holds whenever either $x_{i-1} x_i
x_{i+1}$ is a right angle or $|x_{i-1}-x_i|=|x_i-x_{i+1}|$. In this
"billiard" we
have closed trajectories in the form of arbitrary rectangles inscribed into the
circle. Thus, we may have one pair of sides arbitrarily small.}
\end{remark}
\begin{definition}\label{generic}
Let $X\subset \R^{m+1}$ be a smooth hypersurface.
Then $X$ is {\it $n$-generic} if $L_X: G(X,n)\to \R$ is a Morse function.
\end{definition}

A justification for this definition is provided by the next lemma.
\begin{lemma}\label{gp}
There is a massive subset $E$ in the space of embeddings
$S^m \to \R^{m+1}$
such that for every $f \in E$ the hypersurface Im $f$ is $n$-generic for
all $n$.
\end{lemma}

Recall that the space of smooth maps from one manifold to another is
considered in the Whitney $C^{\infty}$ topology; a { massive} set is a
countable intersection of open
dense sets. Due to the Baire property, a massive set is dense -- see,
\cite{GG}.

{\bf Proof.}
Consider $n$ germs of immersions
$$\phi_i: (S^m,s_i) \to (\R^{m+1},x_i),\quad i=1, \dots, n,$$
and assume that the targets $x_1, \dots, x_n$ satisfy $x_i \neq x_{i+1}$
for $i=1, \dots, n$. Then a germ of the respective perimeter length function
$(S^m,s_1) \times \dots \times (S^m,s_n) \to \R$
is defined:
$$ L(t_1, \dots, t_n) = \sum_1^n |\phi_{i+1}(t_{i+1})- \phi_i (t_i)|.$$
Clearly, the first partial derivatives of $L$ depend on the first
derivatives of  $\phi_i$.
Therefore we may consider the situation on the level of 1-jets. Namely, let
$${\cal U} \subset J^1 (S^m,\R^{m+1}) \times \dots \times J^1
(S^m,\R^{m+1})$$ consist of
multi-jets of immersions satisfying the following three conditions:
\begin{itemize}
\item[{\rm (i)}] the targets $x_1, \dots, x_n$  satisfy $x_i \neq x_{i+1}$
for $i=1,
\dots, n$;
\item[{\rm (ii)}] the vector
\begin{eqnarray}
\nu_i = {{x_{i} - x_{i-1}}\over{|x_{i} - x_{i-1}|}} +
{{x_{i}- x_{i+1}}\over{|x_{i}- x_{i+1}|}} \label {nu}
\end{eqnarray}
 does not vanish for $i=1, \dots, n$;
\item[{\rm (iii)}] the sources $t_1, \dots, t_n$ satisfy $t_i \neq t_j$ for
all $i \neq j$.
\end{itemize}
The first requirement is needed for $L$ to be smooth and the third for the
multi-jet
transversality  theorem to be applicable; the role of the second one will
become clear
shortly. Note that
${\cal U}$ is an open subset of the multi-jet space. Also consider the space of
1-jets of functions $J^1 (S^m
\times \dots \times S^m)$, and let $D$ be its submanifold of codimension
$mn$ consisting  of the
1-jets with trivial differential. Assigning the 1-jet of the
perimeter length function $L$ to a multi-jet in ${\cal U}$ provides a  map
$$\pi: {\cal U} \to J^1 ((S^m)^{\times n}).$$
Claim: {\it the map $\pi$ is transversal to the submanifold $D$}.

Assuming the claim, the proof proceeds as follows. The 1-jet extension of
the perimeter length
function $L$ determines the section $j^1 (L)$ of the 1-jet bundle
$$J^1 ((S^m)^{\times n}) \to (S^m)^{\times n},$$
 and the critical points
of $L$ are non-degenerate if and only if this section is transversal to
$D$. The claim implies that
$\Delta = \pi^{-1} (D)$ is a submanifold in ${\cal U}$. According to the
multi-jet transversality
theorem (see \cite{GG}), there exists a massive set $E_n$ of embeddings
$S^m \to \R^{m+1}$
such that for $f \in E_n$ the
multi-jet  extension $j^1 (f): G(S^m,n) \to  {\cal U}$ is transversal to
$\Delta$.
Since $j^1 (L) = \pi \circ j^1 (f)$, it follows from the claim that
$j^1(L)$ is transversal to $D$
for $f \in E_n$, that is, the critical points of the perimeter length
function are
non-degenerate. Setting $E = \cap_{n=2}^{\infty} E_n$ completes the proof.

It remains to prove the italized claim above.
It is convenient to choose local coordinates in the 1-jet spaces involved.
Consider a multi-jet
$\bar \phi = (\phi_1, \dots, \phi_n) \in \Delta$, and let
$s_i \in S^m$ be the source of $\phi_i$. For each $i$ identify a
neighborhood of $s_i$ with an open
disk $U \subset \R^m$. Then a neighborhood of $\bar \phi$ in
$J^1(S^m,\R^{m+1})^{\times n}$ is
identified with
$J^1 (U, \R^{m+1})^{\times n}$ and $\pi$ becomes a map
$$\pi: J^1 (U, \R^{m+1})^{\times n} \to J^1(U^{\times n}).$$
The 1-jet space $J^1 (U, \R^{m+1})$ consists of triples $(u,x,A)$ where
$u\in U$ is the source,
$x \in \R^{m+1}$ is the target and $A: \R^m \to \R^{m+1}$ is a linear map
(derivative of a map $U\to \R^{m+1})$; the multi-jet space
$J^1 (U, \R^{m+1})^{\times n}$ consists of $n$-tuples
$(u_i,x_i,A_i),\ i=1,\dots, n$, of such
triples. The  space $J^1(U^{\times n})$ consists of 1-jets of functions
$\psi: U^{\times n} \to\R$, that is, of $2n+1$-tuples $(u_i,p_i,z),\ i=1,
\dots, n$,
where $u_i \in U$, $p_i = \partial \psi/\partial u_i \in (\R^{m+1})^*$ is a
covector  and $z =\psi(u_1, \dots, u_n) \in \R$.
In these coordinates we explicitly describe the map $\pi: (u_i,x_i,A_i) \to
(u_i,p_i,z)$:
$$z = \sum_{i=1}^n |x_{i+1} - x_{i}|,\quad p_i (v_i) = \langle \nu_i, A_i
(v_i) \rangle$$
where $v_i \in \R^m$ is a test vector and the vector $\nu_i$ is as in
(\ref{nu}).
 The first  formula is obvious and the second
was established in the proof of Proposition \ref{prop2}. Identifying
vectors and covectors by the
Euclidean structure, one has:
\begin{eqnarray}
p_i = A_i^* (\nu_i). \label{mom}
\end{eqnarray}
Consider a multi-jet $\bar \phi \in \Delta \subset J^1 (U,
\R^{m+1})^{\times n}$; as before,
$\phi_i = (u_i,x_i,A_i)$. We want to show that
\begin{eqnarray}
d\pi (T_{\bar \phi} J^1 (U, \R^{m+1})^{\times n}) + T_{\pi(\bar \phi)} D =
T_{\pi(\bar \phi)} J^1(U^{\times n}). \label{tr}
\end{eqnarray}
The space $T_{\pi(\bar \phi)} D$ consists of vectors whose $p_i$ components
vanish, while the $u_i$
and $z$ components are arbitrary. Thus the equality (\ref{tr}) will follow
once we show that every
vector in
$T_{\pi(\bar \phi)} J^1(U^{\times n})$ with trivial $u_i$ and $z$
components is in the image of $d\pi$.

Consider an infinitesimal deformation $(u_i,x_i,A_i + \varepsilon B_i)$ of
$\bar \phi \in J^1 (U,
\R^{m+1})^{\times n}$ where $B_i : \R^m \to \R^{m+1}$ is a linear map; this
deformation determines
a tangent vector $\xi \in T_{\bar \phi} J^1 (U, \R^{m+1})^{\times n}$.
Since $\pi (\bar \phi) \in
D$, it follows from (\ref{mom}) that $A_i^* (\nu_i) = 0$. Therefore $d\pi
(\xi)$ is a vector
in $T_{\pi(\bar \phi)} J^1(U^{\times n})$ whose $p_i$ component is
$B_i^*(\nu_i)$, while the $u_i$
and $z$ components vanish. Since $\nu_i \neq 0$, the vector
$B_i^*(\nu_i) \in \R^m$ can be made arbitrary by varying $B_i$, and the
result follows.
\proofend

We will use Proposition \ref{prop2} to tackle the problem of finding
topological lower bounds on the number of closed billiard trajectories by
applying the methods of Morse theory. The function $L_X$ being $D_n$-invariant,
we will use  equivariant Morse and Lusternik-Schnirelman theories. Namely, one
has the next result.

\begin{proposition}\label{thm7} Let $X\subset \R^{m+1}$ be a smooth strictly
convex hypersurface.
Then for any odd $n\geq 3$ the number of $D_n$-orbits of $n$-periodic billiard
trajectories in $X$ is greater than the cup-length of
\begin{eqnarray}
H^\ast(G(S^m,n)/D_n;\Z_2).\label{cup}
\end{eqnarray}
If $X$ is $n$-generic then the number of $D_n$-orbits of
$n$-periodic billiard trajectories
in $X$ is not less than the sum of Betti numbers
\begin{eqnarray}
\sum_i \dim_{\Z_2}H^i(G(S^m,n)/D_n;\Z_2).\label{sum}
\end{eqnarray}
\end{proposition}

{\bf Proof.} Start with the following claim: for odd $n$ the cohomology group
of the quotient
$G(S^m,n)/D_n$ coincides with the equivariant cohomology of $G(S^m,n)$:
\begin{eqnarray}
H^j(G(S^m,n)/D_n;\Z_2)\, \simeq \,
H^j(E{D_n}\times_{D_n}G(S^m,n),\Z_2)\label{isom}
\end{eqnarray}
where $E{D_n}$ is a contractible space with a free $D_n$-action.

Indeed, consider the Leray spectral sequence
of the projection $E{D_n}\times_{D_n}G(S^m,n)\to G(S^m,n)/D_n$, see
\cite{Br}. The $E_2$-term
has the form $E_2^{p,q}=H^{p}(G(S^m,n)/D_n;\mathcal L^q)$, where $\mathcal
L^q$ is the Leray sheaf.  Let $\bar x\in G(S^m,n)$ be an orbit. Since $n$
is odd,
no reflection in $D_n$  belongs to the stabilizer $H_{\bar x}$. Thus the
stabilizer of $\bar x$ is a cyclic subgroup  $H_{\bar x}\subset D_n$ of
odd order.  The stalk of $\mathcal L^q$ over the orbit of $\bar x$ is
$$H^q(ED_n\times_{D_n} D_n/H_{\bar x};\Z_2) = H^q(ED_n/H_{\bar x};\Z_2). $$
It follows that $\mathcal L^q=0$ for $q>0$. Therefore the Leray
spectral sequence is nonzero only along the $p$-axis. This implies
(\ref{isom}).

Now we use Proposition \ref{prop2}.
Fix a sufficiently small $\varepsilon>0$ and consider the function $L_X:
G_\varepsilon(X,n)
\to \R$, determined by the billiard hypersurface $X\subset \R^{m+1}$. Since
$L_X$ is
$D_n$-invariant, the set of critical points of $L_X$ is also $D_n$-invariant.
Our task is to estimate the number of critical $D_n$-orbits.
We have to take into account the presence of the boundary.
However, due to statement (d) of Proposition \ref{prop2},
the boundary points make no contribution to the topology of
$G_\varepsilon(X,n)$ -- cf. \cite{BM}. In other words, the principles of the
critical point theory
apply to $L_X: G_\varepsilon(X,n)\to \R$ the same way as if
$G_\varepsilon(X,n)$ were a manifold
without boundary.

Assume that $X$ is $n$-generic. Then, using the negative gradient flow of
$L_X$, we obtain a
$D_n$-equivariant cell decomposition of $G_\varepsilon(X,n)$.  The number
of cells in the resulting
cell decomposition of the quotient space  $G_\varepsilon(X,n)/D_n$ equals
the number of
critical $D_n$-orbits of
$L_X$, i.e., the number of $D_n$-orbits of $n$-periodic billiard
trajectories in $X$.
This proves the second statement of Proposition {\ref{thm7}}.

To prove the first statement,  apply  equivariant
Lusternik - Schnirelman theory developed in \cite{Mar}, see also \cite{CP}.
We  use Theorems 3.2 and 1.13  from \cite{Mar}.
The first example in section 1.14 of \cite{Mar} of singular multiplicative
$D_n$-cohomology theory
is given by $Y\mapsto H^\ast(Y/D_n;\Z_2)$; in our case
$H^\ast(G(X,n)/D_n;\Z_2)$
coincides with the equivariant cohomology $H^\ast_{D_n}(G(X,n);\Z_2)$ due
to (\ref{isom})
and statement (b) of Proposition \ref{prop2}.
\proofend

\section{Equivariant cohomology of the cyclic configuration space of the
sphere}

In view of Proposition \ref{thm7}, in order to find lower bounds on the
number of
periodic trajectories of billiards we need to compute the equivariant
cohomology ring of the cyclic configuration space $G(S^m,n)$.
This is the main goal of this section.

Our computation of the equivariant cohomology will be based on Theorem
\ref{thm3}, giving the structure of the usual cohomology ring. We will also
use the
method which was
first suggested by I. Babenko \cite{Ba}. It consists in applying Morse
theory in the opposite direction, that is, studying the topology of the cyclic
configuration
space $G(S^m,n)$ by examining the billiard inside the round ball in
$\R^{m+1}$; in
the latter case the closed trajectories are readily described. Note that a
similar idea of studying the closed geodesics on the round sphere leads to a
computation of the homology of the loop space of the sphere -- see
\cite{Bo} and the references
therein.

{}From now on we will assume that the number $n$ is odd and that $m>1$.

Let $X = S^m \subset {\R}^{m+1}$ be the unit sphere. Consider the corresponding
length function $L_X: G(X,n)\to \R$. The critical points of $L_X$
are precisely the closed
billiard trajectories inside $X$ having $n$ reflections.
Each such $n$-periodic trajectory lies in a two-dimensional plane $P$
 through the center of the sphere $X$. The intersection $P\cap X$
is a unit circle and the reflections go  the same way as in the
plane circular billiard $P\cap X$. Hence  any billiard
trajectory $(x_1, x_2, \dots, x_n)$ is a plane regular $n$-gon,
possibly star-shaped, inscribed into the circle $P\cap X$. If $n$ is not
prime then such a polygon may be multiple, i.e., it may traverse itself
several times.
The angle $\alpha$ between $x_1$ and $x_2$ is of the form $\alpha=2\pi r/n$,
where $1\le r\le (n-1)/2$.
This number is clearly related to the rotation number (cf. Example \ref{ex1}).
The following picture shows the $7$-periodic trajectories:
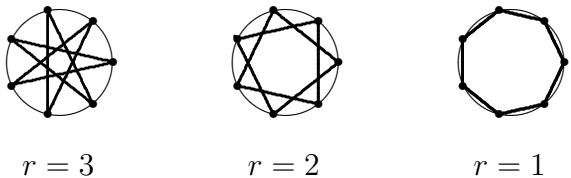
\begin{figure}[h]
 \begin{center}
\setlength{\unitlength}{1cm}
\begin{picture}(9,4)
\put(2,2){\circle{3}}
\put(2.4429, 2.5548){\circle*{0.1}}
\put(1.8426,2.6923){\circle*{0.1}}
\put(1.3608, 2.3089){\circle*{0.1}}
\put(2.71, 2){\circle*{0.1}}
\put(2.4429, 1.4452){\circle*{0.1}}
\put(1.8426,1.3077){\circle*{0.1}}
\put(1.3608, 1.6911){\circle*{0.1}}
\multiput(1.8426,2.6923)(0.006,-0.0124){100}{\circle*{0.01}}
\multiput(2.4429, 1.4452)(-0.0108,0.0086){100}{\circle*{0.01}}
\multiput(1.3608, 2.3089)(0.0134,-0.0030){100}{\circle*{0.01}}
\multiput(2.71, 2)(-0.0135,-0.0030){100}{\circle*{0.01}}
\multiput(1.3608, 1.6911)(0.0108,0.0086){100}{\circle*{0.01}}
\multiput(2.4429, 2.5548)(-0.0060,-0.0124){100}{\circle*{0.01}}
\multiput(1.8426,1.3077)(0.00,0.0133){100}{\circle*{0.01}}


\put(5,2){\circle{3}}
\put(5.4429, 2.5548){\circle*{0.1}}
\put(4.8426,2.6923){\circle*{0.1}}
\put(4.3608, 2.3089){\circle*{0.1}}
\put(5.71, 2){\circle*{0.1}}
\put(5.4429, 1.4452){\circle*{0.1}}
\put(4.8426,1.3077){\circle*{0.1}}
\put(4.3608, 1.6911){\circle*{0.1}}
\multiput(5.71, 2)(-0.00867,0.006923){100}{\circle*{0.01}}
\multiput(5.71, 2)(-0.00867,-0.006923){100}{\circle*{0.01}}
\multiput(4.8426,1.3077)(-0.005182,0.01053){100}{\circle*{0.01}}
\multiput(4.3608, 2.3089)(0.01082,0.002459){100}{\circle*{0.01}}
\multiput(5.4429, 2.5548)(0,-0.0110){100}{\circle*{0.01}}
\multiput(5.4429, 1.4452)(-0.0108,0.00245){100}{\circle*{0.01}}
\multiput(4.3608, 1.6911)(0.004818,0.01){100}{\circle*{0.01}}

\put(8,2){\circle{3}}
\put(8.4429, 2.5548){\circle*{0.1}}
\put(7.8426,2.6923){\circle*{0.1}}
\put(7.3608, 2.3089){\circle*{0.1}}
\put(8.71, 2){\circle*{0.1}}
\put(8.4429, 1.4452){\circle*{0.1}}
\put(7.8426,1.3077){\circle*{0.1}}
\put(7.3608, 1.6911){\circle*{0.1}}
\multiput(8.71, 2)(-0.00267,0.005548){100}{\circle*{0.01}}
\multiput(8.71, 2)(-0.00267,-0.005548){100}{\circle*{0.01}}
\multiput(8.4429, 2.5548)(-0.006,0.00137){100}{\circle*{0.01}}
\multiput(8.4429, 1.4452)(-0.006,-0.00137){100}{\circle*{0.01}}
\multiput(7.8426,1.3077)(-0.00481,0.003834){100}{\circle*{0.01}}
\multiput(7.8426,2.6923)(-0.00481,-0.003834){100}{\circle*{0.01}}
\multiput(7.3608, 1.6911)(0,0.006178){100}{\circle*{0.01}}

\put(1.5,0.5){$r=3$}
\put(4.5,0.5){$r=2$}
\put(7.5,0.5){$r=1$}

\end{picture}
 \end{center}
\caption{Critical submanifolds for $G(S^n,7)$}
\end{figure}

Note that any two $n$-periodic trajectories in $X$ with the same rotation
number $r$ can be continuously deformed one to another. We conclude that
the critical
points of the function
$L_X: G(X,n)\to \R$ form a disjoint union of connected submanifolds
\begin{eqnarray}
V_0, V_1, \dots, V_{(n-3)/2},\label{rot}
\end{eqnarray}
where $V_p$ denotes the set of all closed trajectories with the rotation
number $(n-1-2p)/2$.
Each $V_p$ is diffeomorphic to the Stiefel manifold $V_{2, m+1}$ and hence
is a closed manifold of dimension $2m-1$. The next result is due to Babenko
\cite{Ba}.

\begin{proposition}\label{propbab}
 If $X = S^m\subset \R^{m+1}$ is a round sphere then
\begin{enumerate}
\item[(a)] The function
$L_X: G(X,n)\to \R$ is nondegenerate in the sense of Bott.
\item[(b)] The index of the critical manifold
$V_p$ equals $2p(m-1)$.
\item[(c)] The critical values of the function $L_X$ on the critical
manifolds $V_p$
increase: $L_X(V_p)<L_X(V_{p'})$ for $p<p'$.
\end{enumerate}
\end{proposition}

As an important addition to Proposition \ref{propbab} we make the following
observation.

\begin{proposition}\label{perf}
If $X$ is a round sphere then  $L_X: G(X, n)\to \R$
is a perfect Bott function with respect to the field $\Z_2$.
\end{proposition}

{\bf Proof.} Let us first explain the meaning of our statement.
Choose $\varepsilon >0$ as in
Proposition \ref{prop2} holds. We may find constants $c_0, \dots,
c_{(n-3)/2}$ with
$$L_X(V_p)<c_p<L_X(V_{p+1}).$$
Set: $F_p=L_X^{-1}((-\infty, c_p])$. We obtain a filtration
\begin{eqnarray}
F_0\subset F_1\subset \dots \subset F_{(n-3)/2}=
G_\varepsilon(X,n),\label{filtration}
\end{eqnarray}
and our statement means that the sum of the Poincar\'e polynomials of the
pairs $(F_p, F_{p-1})$ with $\Z_2$ coefficients equals
the Poincar\'e polynomial of $G_\varepsilon(X,n)$.

Indeed, the mod 2 Poincar\'e polynomial of $(F_p, F_{p-1})$  equals
\begin{eqnarray}
\sum_{j} t^j\dim_{\Z_2}H^j(F_p, F_{p-1};\Z_2) \, =\,
t^{2p(m-1)}[t^{2m-1}+t^m +t^{m-1} +1].\label{poinc}
\end{eqnarray}
Here we used the fact that $(F_p, F_{p-1})$ is homotopy equivalent to the
Thom space
of the negative normal bundle of $V_p$, the Thom isomorphism, the index
computation, given by
Proposition \ref{propbab}, and  the fact that the  Poincar\'e polynomial
with coefficients in
$\Z_2$  of the Stiefel manifold
$V_{2, m+1}$ is $t^{2m-1}+t^m +t^{m-1} +1$, cf. \cite{Br1}. Summing
formulae (\ref{poinc})
for all $p=0, \dots, (n-3)/2$ we obtain
$$
\begin{array}{l}
\displaystyle
\sum_{p=0}^{(n-3)/2}\sum_{j} t^j\dim_{\Z_2} H^j(F_p, F_{p-1};\Z_2) \, =\, \\ \\
\displaystyle
=\, (t^{2m-1}+t^m +t^{m-1}+1)\cdot \frac{t^{(n-1)(m-1)} -1}{t^{2(m-1)}-1}
\, =\, \\ \\
\displaystyle
=\, \frac{(t^m+1)(t^{(n-1)(m-1)}-1)}{t^{(m-1)}-1},\\
\end{array}
$$
which, according to Theorem \ref{thm3}, coincides with the Poincar\'e
polynomial of $G(X,n)$.
\proofend

Propositions \ref{propbab} and \ref{perf} hold for even $n$ as well, but
then their statements
are slightly different. We will not need these results in this paper.

\begin{remark}\label{rmk8} {\rm For $m>3$ there exists
a different proof of Proposition \ref{perf}, which does not use Theorem
\ref{thm3} and provides an independent computation of the Poincar\'e
polynomial of the
cyclic configuration space
of the sphere. It is quite straightforward. One considers the spectral
sequence of the filtration  (\ref{filtration}), where
$$E_1^{p,q}=H^{p+q}(F_p, F_{p-1};\Z_2).$$
The calculation of the relative homology
$H^{p+q}(F_p, F_{p-1};\Z_2)$ as in the above proof
and elementary geometric considerations show that all differentials must be
zero provided $m>3$.

This argument fails  for $m=2$ and for $m=3$.
The nonzero terms and the differential $d_1$ of the spectral sequence for
$m=2$ and $m=3$ are shown
in Figure 4. Unlike the case $m>3$, this picture does not imply that the
differentials vanish.

Note also that for $n=3$ the proof of Proposition \ref{perf} trivializes:
the filtration
(\ref{filtration}) consists of a single term only.}\end{remark}

\setlength{\unitlength}{0.9cm}
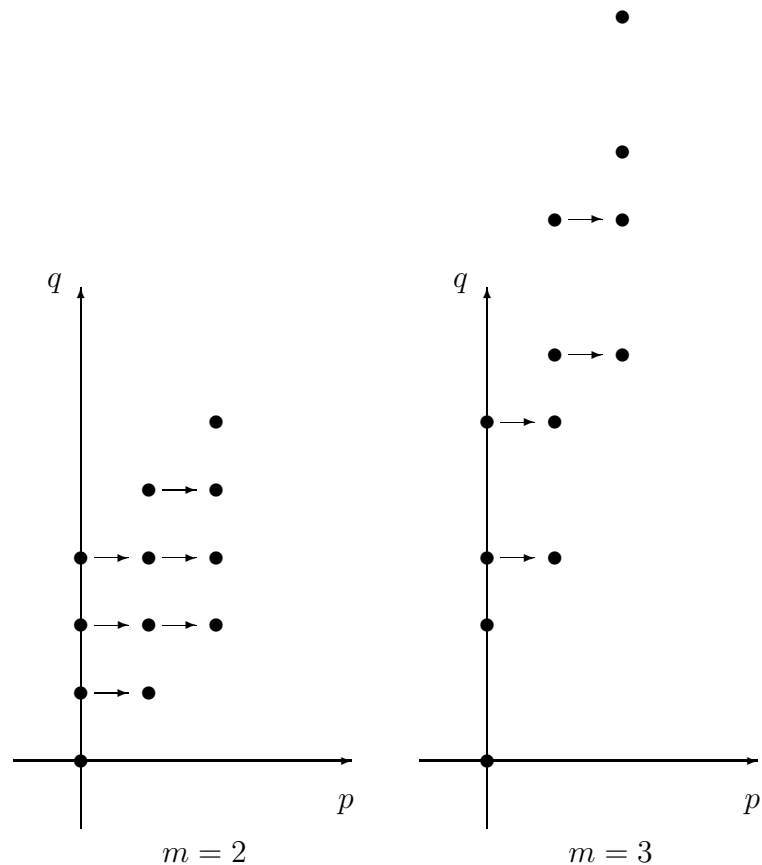
\begin{figure}[h]
 \begin{center}
\begin{picture}(10, 12)
\linethickness{0.15mm}
\put(0,1){\vector(1,0){5}}
\put(1,0){\vector(0,1){8}}
\put(1,1){\circle*{0.2}}
\put(1,2){\circle*{0.2}}
\put(1,3){\circle*{0.2}}
\put(1,4){\circle*{0.2}}
\put(2,3){\circle*{0.2}}
\put(2,4){\circle*{0.2}}
\put(2,5){\circle*{0.2}}
\put(2,2){\circle*{0.2}}
\put(3,3){\circle*{0.2}}
\put(3,4){\circle*{0.2}}
\put(3,5){\circle*{0.2}}
\put(3,6){\circle*{0.2}}
\put(4.8,0.3){$p$}
\put(0.5,8){$q$}
\linethickness{0.05mm}
\put(1.2,3){\vector(1,0){0.5}}
\put(1.2,4){\vector(1,0){0.5}}
\put(1.2,2){\vector(1,0){0.5}}
\put(2.2,3){\vector(1,0){0.5}}
\put(2.2,4){\vector(1,0){0.5}}
\put(2.2,5){\vector(1,0){0.5}}
\linethickness{0.15mm}
\put(6,1){\vector(1,0){5}}
\put(7,0){\vector(0,1){8}}
\put(7,1){\circle*{0.2}}
\put(7,3){\circle*{0.2}}
\put(7,4){\circle*{0.2}}
\put(7,6){\circle*{0.2}}
\put(8,4){\circle*{0.2}}
\put(8,6){\circle*{0.2}}
\put(8,7){\circle*{0.2}}
\put(8,9){\circle*{0.2}}
\put(9,7){\circle*{0.2}}
\put(9,9){\circle*{0.2}}
\put(9,10){\circle*{0.2}}
\put(9,12){\circle*{0.2}}
\put(10.8,0.3){$p$}
\put(6.5,8){$q$}
\linethickness{0.05mm}
\put(7.2,4){\vector(1,0){0.5}}
\put(7.2,6){\vector(1,0){0.5}}
\put(8.2,7){\vector(1,0){0.5}}
\put(8.2,9){\vector(1,0){0.5}}
\put(2.2,-0.5){$m=2$}
\put(8.2,-0.5){$m=3$}
\end{picture}
\end{center}
\caption{Nonzero terms of the spectral sequence for $m=2$ and $m=3$.}
\end{figure}

The main result of this section is the next theorem.

\begin{theorem}\label{thm61}
Let $m \geq 3$ and let $n$ be an odd integer. Then the cohomology ring
$H^\ast(G(S^m;n)/D_n;\Z_2)$
is multiplicatively generated by cohomology classes
$$\sigma_{2i}\in H^{2i(m-1)}(G(S^m;n)/D_n;\Z_2), \quad where \quad i=1, \,
2, \dots$$
and by classes
$$e\in H^1(G(S^m;n)/D_n;\Z_2)\quad and \quad u\in H^m(G(S^m;n)/D_n;\Z_2).$$
These classes satisfy the following relations
\begin{itemize}
\item[{\rm (i)}] $\sigma_i=0$ for $i\ge n-1$;
\item[{\rm (ii)}] $\sigma_i\sigma_j = \displaystyle{i+j \choose
i}\sigma_{i+j}+ \epsilon_{ij}\sigma_{i+j-2}ue^{m-2}$,
where $i$ and $j$ are even and $\epsilon_{ij}\in \Z_2$;
\item[{\rm (iii)}] $e^m=0$;
\item[{\rm (iv)}] $u^2=0$.
\end{itemize}
The Poincar\'e polynomial of the quotient space $G(S^m,n)/D_n$ with
coefficients in $\Z_2$ equals
$${{(t^{(n-1)(m-1)}-1)}\over{(t^{2(m-1)}-1)}}\cdot {{t^m-1}\over{t-1}}\cdot
(t^m+1),$$
and the sum of Betti numbers is $m(n-1)$.
\end{theorem}

The rest of this section consists of the proof of Theorem \ref{thm61}.

Start with the following lemma.

\begin{lemma}\label{lm3}
Let $G$ be a finite group acting simplicially on a finite polyhedron $Y$
such that the action of $G$
on the cohomology $H^\ast(Y;\Z_2)$ is trivial. Suppose that $G'\subset G$
is a subgroup of odd index such that $G$ acts trivially on $H^\ast(G';\Z_2)$.
Then the induced homomorphism of the equivariant cohomology
\begin{eqnarray}
H^\ast_G(Y;\Z_2) \to H^\ast_{G'}(Y;\Z_2)\label{isom1}
\end{eqnarray}
is an isomorphism.
\end{lemma}
The notation $H^\ast_G(Y;\Z_2)$ stands for the equivariant cohomology
$H^\ast(EG\times_G Y;\Z_2)$.

{\bf Proof.} We will use the comparison theorem for spectral sequences.
Since $G$ acts trivially on $H^\ast(Y;\Z_2)$, the Serre
spectral sequence of the
fibration $EG\times_G Y\to BG$ with fiber $Y$, converging to
$H^\ast_G(Y;\Z_2)$, has the initial term
$$E_2^{p,q} \, =\, H^p(G;\Z_2)\otimes H^q(Y;\Z_2).$$
Similarly, we have a spectral sequence with the initial term
$${E'}_2^{p,q} \, =\, H^p(G';\Z_2)\otimes H^q(Y;\Z_2)$$
converging to  $H^\ast_{G'}(Y;\Z_2)$. The inclusion $G'\to G$ induces
a homomorphism of the spectral sequences $E\to E'$ which is
an isomorphism of the $E_2$-terms (cf. \cite{Bw}, Proposition 10.4, chapter
3).
Hence, by the comparison theorem for spectral sequences,
 (\ref{isom1}) is an isomorphism. \proofend

Next we compute the $D_n$-equivariant cohomology of the critical manifolds
of the function $L_X$.

\begin{proposition}\label{rot1} Suppose that $m \geq 3$ and $n$ is odd.
Let $V_p$ be any of the critical submanifolds (\ref{rot}) with the induced
action of the dihedral
group $D_n$. Then
\begin{itemize}
\item[(a)]
the equivariant cohomology
ring $H^\ast_{D_n}(V_p;\Z_2)$ has two multiplicative generators
$$e\in H^1_{D_n}(V_p;\Z_2)\quad and \quad u\in H^m_{D_n}(V_p;\Z_2);$$
\item[(b)] they satisfy the relations
$e^m=0$ and $u^2=0;$
\item[(c)]
the classes $e^iu^j$, with $i=0, 1, \dots, m-1$ and $j=0, 1$ form an
additive basis of $H^\ast_{D_n}(V_p;\Z_2)$;
\item[(d)] the canonical homomorphism
$H^\ast(V_p/D_n;\Z_2) \to H^\ast_{D_n}(V_p;\Z_2)$ is an isomorphism;
\item[(e)] the kernel of the canonical homomorphism
$\Phi: H^\ast(V_p/D_n;\Z_2)\to H^\ast(V_p;\Z_2)$
coincides with the ideal generated by $e$, i.e.,  is the
linear span of the classes $e^iu^j$ where $i\ge 1$.
\end{itemize}
\end{proposition}

{\bf Proof.} First note that (d) follows from the isomorphism (\ref{isom})
in the proof of Proposition
 \ref {thm7}. Indeed, the stabilizer of any orbit in  $V_p$ is a cyclic
subgroup of odd order.

To compute the equivariant cohomology ring $H^\ast_{D_n}(V_p;\Z_2)$ we will
apply Lemma
\ref{lm3} with $G=D_n$, $Y=V_p$ and  $G'\simeq \Z_2$ being a subgroup of
the dihedral
 group $D_n$ generated by a reflection. Since $V_p$ is homeomorphic to the
Stiefel manifold $V_{2,
m+1}$, the cohomology of $V_p$ with $\Z_2$ coefficients is isomorphic to
$\Z_2$ in dimensions $0, m-1, m, 2m-1$ and  is trivial in all other
dimensions.
It follows that $D_n$ acts
trivially on $H^\ast(V_p;\Z_2)$; likewise, $D_n$ acts trivially on
$H^\ast(G';\Z_2)$.
Applying Lemma \ref{lm3}, we conclude that there exists a ring isomorphism
\begin{eqnarray}
H^\ast_{D_n}(V_p;\Z_2)\simeq H^\ast_{G'}(V_p;\Z_2)\simeq
H^\ast(V_p/G';\Z_2).\label{factor}
\end{eqnarray}
The second isomorphism follows since the reflection acts freely on $V_p$
for odd $n$. Formulae
(\ref{factor}) show that, computing the equivariant cohomology
$H^\ast_{D_n}(V_p;\Z_2)$, we may ignore a large part of the $D_n$-action
and keep track  only of the $G'$-action.

Note that $V_p$ can be identified with
the variety of ordered pairs $(v_1, v_2)$ of unit vectors in $\R^{m+1}$
making the
angle of $\alpha_p=\pi(n-1-2p)/n$. The reflection (i.e. the generator of
$G'$) acts on such pairs by
sending $(v_1, v_2)$ to $(v_1, v'_2)$, where
$v'_2$ is the reflection of $v_2$ in the line spanned by $v_1$.

Applying the Gram - Schmidt orthogonalization,
construct a diffeomorphism between the quotient space $V_p/G'$ and the set
of pairs $(v_1, v_2)$ of
mutually orthogonal unit vectors in $\R^{m+1}$ with identification $(v_1,
v_2) \simeq (v_1, -v_2)$. In
other words, $V_p/G'$ is diffeomorphic to space of pairs $(v, \ell)$, where
$e\in \R^{m+1}$ is a unit
vector, and  $\ell\subset \R^{m+1}$ is a one-dimensional linear subspace
orthogonal to $v$.
The projection $(v,\ell)\mapsto v$ identifies $V_p/G'$ with the
projective tangent bundle of $S^m$.

Alternatively, projecting $(v,\ell)\mapsto \ell$, we  view  $V_p/G'$ as the
space of a
unit sphere bundle of a rank $m$ vector bundle $\xi$ over the projective
space $\RP^m$. The fiber of $\xi$
over a line $\ell\in \RP^m$ is the orthogonal complement $\ell^\perp$ of
$\ell$. The spectral sequence
of this unit sphere bundle
$$E_2^{p\, q}\, =\, H^p(\RP^m;\Z_2)\otimes H^q(S^{m-1};\Z_2) \Rightarrow
H^{p+q}(V_p/G';\Z_2)$$
has only two rows and the only possibly nontrivial differential is the
transgression
$$d_m: E_m^{0, m-1}=H^{m-1}(S^{m-1};\Z_2)\to E_m^{m,0}=H^m(\RP^m;\Z_2).$$

We claim that the differential $d_m: E_m^{0, m-1}\to E_m^{m,0}$ is an
isomorphism. The image
of the generator of $H^{m-1}(S^{m-1};\Z_2)$ under $d_m$ is the top
Stiefel-Whitney class
$w_m(\xi)\in H^m(\RP^m;\Z_2)$. Let $\eta$ be the tautological line bundle
over $\RP^m$
whose fiber over a line $\ell$ is $\ell$ itself. Then the Whitney
sum $\xi\oplus \eta$ is the trivial bundle of rank $m+1$. Since the total
Stiefel-Whitney class of $\eta$
is $1+e$,  where $e\in H^1(\RP^m;\Z_2)$ is the generator, the Cartan's formula
\begin{eqnarray}
(1+e)(1 +w_1(\xi) +\dots +w_m(\xi)) = 1
\end{eqnarray}
gives $w_j(\xi) = e^j,$ for all $ j=1, \dots, m.$  In particular,
$w_m(\xi)=e^m$.

These arguments  completely describe
the ring structure of $H^\ast(V_p/G';\Z_2)\simeq H^\ast_{D_n}(V_p;\Z_2)$.
Namely,
in the above spectral sequence, the classes $1, e, \dots, e^{m-1}$ survive
in the bottom row;
also there is a class $u\in E^{1,m-1}_{\infty}$
such that the nonzero classes in row $q=m-1$, surviving in $E_\infty$, are
$u, ue, \dots, ue^{m-1}$.
The cohomology classes $e$ and $u$ lift uniquely from $E_\infty$ to
$H^\ast(V_p/G';\Z_2)$ and satisfy the
same relations therein.  This proves statements (a), (b), (c).

It remains to prove statement (e). Clearly, $\Phi$ is a ring homomorphism
and $\Phi(e)=0$. Therefore
(e) will follow once we show that $\Phi: H^m(V_p/G';\Z_2) \to
H^m(V_p;\Z_2)$ is an isomorphism in degree
$m$.

Consider the product $V_p\times [0,1]$ and identify the points $(x,1)$ and
$(Tx,1)$ for all $x\in V_p$,
where $T\in G'$ denotes the reflection. After the identification, we obtain
a compact $2m$-dimensional
manifold with boundary $Y$, so that $Y$ is homotopy equivalent to $V_p/G'$
and $\partial Y$ is
diffeomorphic to $V_p$. The restriction homomorphism $H^j(Y;\Z_2)\to
H^j(\partial Y;\Z_2)$ coincides with
$\Phi$. Using the Poincar\'e duality we obtain
\begin{eqnarray}
H^j(Y,\partial Y;\Z_2) \simeq H^{2m-j}(Y;\Z_2) \simeq \Z_2\quad {\rm for\
all}\ j.\label{poincare}
\end{eqnarray}
Consider the exact cohomological sequence of $(Y,\partial Y)$. The
homomorphisms $H^{j}(Y;\Z_2) \to
H^{j}(\partial Y;\Z_2)$ are zero for $j=m-1$ and $j=m+1$, therefore the
homomorphisms
$\Z_2=H^{m-1}(\partial Y;\Z_2) \to H^{m}(Y, \partial Y;\Z_2)=\Z_2$ and
$\Z_2=H^{m+1}(Y, \partial Y;\Z_2)\to
H^{m+1}(Y;\Z_2)=\Z_2$ are isomorphisms. It follows that $\Phi$ is an
isomorphism
 in the following exact sequence (we suppress the coefficients $\Z_2$ from
notation):
$$H^{m-1}(\partial Y) \stackrel \simeq \to H^{m}(Y,\partial Y)\to
H^{m}(Y)\stackrel \Phi\to
H^{m}(\partial Y)\to H^{m+1}(Y, \partial Y)\stackrel \simeq \to H^{m+1}(Y).$$
The proposition is proved.
\proofend

{\bf Proof of Theorem \ref{thm61}.} Let $X\subset \R^{m+1}$ be a round sphere.
Consider filtration (\ref{filtration}) of the space
$G_\varepsilon(X,n)$, where $\varepsilon>0$ is small enough, so that the
claims of
Proposition \ref{prop2} hold.
Denote the space $F_i/D_n$ by $F'_i$. Hence we obtain a filtration
\begin{eqnarray}
F'_0\subset F'_1\subset \dots \subset F'_{(n-3)/2}=
G_\varepsilon(X,n)/D_n\label{filtration1}
\end{eqnarray}
of the space of $D_n$-orbits.

Formulate  the following {\it inductive hypothesis} $\mathcal F_p$,
depending on a number $p=0, 1, \dots, (n-3)/2$:

{\it The cohomology ring $H^\ast(F'_p;\Z_2)$
is multiplicatively generated by cohomology classes
$$\sigma_{2i}\in H^{2i(m-1)}(F'_p;\Z_2), \quad where \quad i=1, \, 2, \dots$$
and by  classes
$$e\in H^1(F'_p;\Z_2)\quad and \quad u\in H^m(F'_p;\Z_2).$$
These classes
satisfy the following relations
\begin{itemize}
\item[{\rm (i)}] $\sigma_i=0$ for $i > 2p$;
\item[{\rm (ii)}] $\sigma_i\sigma_j =
\displaystyle{i+j \choose i}\sigma_{i+j}+\epsilon_{ij}\sigma_{i+j-2}ue^{m-2}$,
where $i$ and $j$ are even and $\epsilon_{ij}\in \Z_2$;
\item[{\rm (iii)}] $e^m=0$;
\item[{\rm (iv)}] $u^2=0$.
\end{itemize}
The Poincar\'e polynomial of the quotient space $F'_p$ with coefficients in
$\Z_2$ equals
$${{(t^{2(p+1)(m-1)}-1)}\over{(t^{2(m-1)}-1)}}\cdot
{{t^m-1}\over{t-1}}\cdot (t^m+1),$$
and the sum of Betti numbers is $2m(p+1)$.

The kernel of the canonical homomorphism
$\Phi: H^\ast(F'_p;\Z_2)\to H^\ast(F_p;\Z_2)$
coincides with the ideal generated by $e$.}

Our aim is to show that statement $\mathcal F_p$ holds for $p=(n-3)/2$.
This would imply Theorem \ref{thm61}.
Argue by induction. Proposition \ref{rot1} implies that $\mathcal F_0$
holds. Hence we need to show that $\mathcal F_p$ implies $\mathcal F_{p+1}$.

Assuming that $\mathcal F_p$ is satisfied, consider the boundary homomorphism
\begin{eqnarray}
\delta : H^i(F'_p;\Z_2) \to H^{i+1}(F'_{p+1}, F'_p;\Z_2).\label{boundary}
\end{eqnarray}
We claim that this homomorphism is trivial for all $i$.
Indeed, the group $H^i(F'_p;\Z_2) $ is nonzero only for
$i\le 2(p+1)(m-1)+1$ (by the assumption $\mathcal F_p$)
and the group $H^{i+1}(F'_{p+1}, F'_p;\Z_2)$ is nonzero only for
$i+1\ge 2(p+1)(m-1)$. The latter follows since $F'_{p+1}/F'_p$ is homotopy
equivalent to the Thom space of a vector
bundle of rank $2(p+1)(m-1)$ over $V_{p+1}/D_n$,
the space of $D_n$-orbits of the critical manifold $V_{p+1}$.
Therefore homomorphism (\ref{boundary})
can be nonzero only for three values of $i$, namely for
$$i=2(p+1)(m-1)-1,\
i=2(p+1)(m-1) \quad {\rm and} \quad i=2(p+1)(m-1)+1.$$

Let us first show that $\delta$ vanishes for $i=2(p+1)(m-1)-1$, i.e., in the
lowest possible
dimension. Consider the critical manifold $V_{p+1}\subset F_{p+1}$ and the
decomposition of its normal bundle
into the negative and positive parts with respect to the Hessian of
function $L_X$. This decomposition
is $D_n$-equivariant, and after factorization by the action of the dihedral
group $D_n$ it
produces two bundles (the negative and the positive) over $V_{p+1}/D_n$.
The space $F'_{p+1}/F'_p$
is homotopy equivalent to the Thom space of the negative bundle.
Consider the commutative diagram
$$
\begin{array}{ccc}
\displaystyle
H^{i}(F'_p;\Z_2) & \stackrel \delta \longrightarrow &
H^{i+1}(F'_{p+1}, F'_{p}; \Z_2)\\ \\
\Phi \downarrow & & \downarrow \Phi_1\\ \\
H^{i}(F_p;\Z_2) & \stackrel {\delta_1} \longrightarrow &
H^{i+1}(F_{p+1}, F_{p}; \Z_2).
\end{array}
$$
Here $\Phi_1$ is defined similarly to $\Phi$, i.e. it is induced by the
canonical projection
$(F_{p+1}, F_p)\to (F_{p+1}',F'_p).$
We know that $\delta_1$ vanishes (by Proposition \ref{perf}). Hence, to
show that $\delta$ vanishes for
$i=2(p+1)(m-1)-1$, it suffices to show that $\Phi_1$ is a monomorphism for
this value of $i$.
We have the following commutative diagram with the vertical maps being the
Thom isomomorphisms
$$
\begin{array}{ccc}
\displaystyle
H^{i+1}(F'_{p+1},F'_p;\Z_2) & \stackrel {\Phi_1} \longrightarrow &
H^{i+1}(F_{p+1}, F_{p}; \Z_2)\\ \\
\simeq \quad\downarrow & & \downarrow \quad \simeq\\ \\
H^{0}(V'_{p+1};\Z_2) & \stackrel {\Phi_2} \longrightarrow &
H^{0}(V_{p+1}; \Z_2).
\end{array}
$$
It is clear that the homomorphism $\Phi_2$ is an isomorphism
(compare statement (e) of Proposition \ref{rot1}), and therefore $\Phi_1$
is also an isomorphism for
$i=2(p+1)(m-1)-1$.

Hence  $\delta$ vanishes for
$i=2(p+1)(m-1)-1$, in other words,
$\delta(\sigma_{2p}ue^{m-3})=0$. Recall that we assume that $m\ge 3$.
Now we will show that $\delta$ vanishes in the two other dimensions as well:
$$\delta(\sigma_{2p}ue^{m-2})=0, \quad \delta(\sigma_{2p}ue^{m-1})=0.$$
It is clear that there exists a class $\tilde e\in H^1(F'_{p+1};\Z_2)$ such
that
$\displaystyle \tilde e|_{F'_p}=e$. Using statement 12 in \S 6, chapter 5
of \cite{Sp}, we obtain
$$\delta(\sigma_{2p}ue^{m-2}) = \delta(\sigma_{2p}ue^{m-3})\cdot \tilde e
=0.$$
Similarly, $\delta(\sigma_{2p}ue^{m-1}) = 0$.

The vanishing of the boundary homomorphism (\ref{boundary}) means that
filtration (\ref{filtration1})
is also perfect, that is, we have an isomorphism
\begin{eqnarray}
H^\ast(F'_{p+1};\Z_2) \, \simeq \, H^\ast(F'_p;\Z_2)\oplus
H^\ast(F'_{p+1}/F'_p;\Z_2)\label{isom2}
\end{eqnarray}
The additive structure of the relative homology $H^\ast(F'_{p+1}/F'_p;\Z_2)$
is given by Proposition \ref{rot1} with a shift of
all degrees by $2(p+1)(m-1)$. Here we use the Thom isomorphism and
the equality between the equivariant cohomology and
the cohomology of the factor space $V_{p+1}/D_n$, which holds by statement
(d) of
Proposition \ref{rot1}.
Hence (\ref{isom2}) fully describes the additive structure of the cohomology
$H^\ast(F'_{p+1};\Z_2)$, which coincides with the statement of the hypothesis
$\mathcal F_{p+1}$.

Now we want to show that the multiplicative structure of
$H^\ast(F'_{p+1};\Z_2)$ is as stated in the
hypothesis $\mathcal F_{p+1}$. Notice first that the cohomology classes
$e$, $u$, $\sigma_2, \dots, \sigma_{2p}$ in $H^\ast(F'_p;\Z_2)$ extend
uniquely to
cohomology classes
in $H^\ast(F'_{p+1};\Z_2)$ and  denote the extensions by the same
symbols. Our next
problem is to identify the class
$\sigma_{2p+2}\in H^{2(p+1)(m-1)}(F'_{p+1};\Z_2)$; caution must be exercised
since the cohomology of $F'_{p+1}$ in this degree is two-dimensional.

Consider again the critical manifold $V_{p+1}\subset F_{p+1}$ and the
$D_n$-equivariant
decomposition of its normal bundle
into the negative and positive parts with respect to the Hessian of the
function $L_X$. After factorization by the action of the dihedral
group, these bundles give two bundles (the negative and the positive) over
$V_{p+1}/D_n$;
the space $F'_{p+1}/F'_p$
is homotopy equivalent to the Thom space of the negative bundle. The Thom
class of the negative
bundle lies in $H^{2(p+1)(m-1)}(F'_{p+1}/F'_p;\Z_2)$, but the last group is
canonically embedded into
$H^{2(p+1)(m-1)}(F'_{p+1};\Z_2)$ via (\ref{isom2}). Hence we define the
class $\sigma_{2p+2}$ as representing this Thom class.

All the generators having been
defined, we want to check that the hypothesis $\mathcal F_{p+1}$
is satisfied. First we note that $\sigma_{2p+2}u$ and $\sigma_{2p+2}e^{m-1}$
are nonzero cohomology classes in $H^\ast(F'_{p+1};\Z_2)$. This would follow
from the Thom Isomorphism Theorem once we show that $e^{m-1}$ and $u$
restrict to nontrivial cohomology classes on $V_{p+1}/D_n$.
Nontriviality of the restriction of $e^{m-1}$ is almost obvious;
indeed, for any $(m-2)$-connected $D_n$-invariant
subset $A\subset F_{p+1}$ we have $e^{m-1}|_{A/D_n}\ne 0$, as
follows by considering the Serre spectral sequence.
In order to show that $u|_{V_{p+1}/D_n}\ne 0$, it is enough to show
that $u|_{V_{p+1}}\ne 0$ (by statement (e) of Proposition \ref{rot1}). Let
$W=V_{3, m+1}$
be the Stiefel manifold
 of triples $(e_1, e_2, e_3)$ of mutually orthogonal unit vectors in
$\R^{m+1}$.
For any $p=0, 1, \dots, (n-3)/2$ denote by $I_p: W\to G(S^m,n)$ the
following map
$$I_p(e_1, e_2, e_3) = (x_1, x_2, \dots, x_n),\quad x_{j+1} =
\cos(j\alpha_p)e_1 +
\sin(j\alpha_p)e_2, $$
where $j=0, \dots, n-1$ and $\alpha_p=\pi(n-1-2p)/n$.
The image of $I_p$ is the critical submanifold $V_p$.
Hence, it is enough to show that $I^\ast_{p}(u)\ne 0$ for any $p$.
Construct a homotopy between $I_p$ and $I_0$. Let
$$H_\tau: W\to G(S^m,n), \quad \tau \in [0,1],$$
be defined by
$H_\tau(e_1, e_2, e_3) = (x_1, x_2, \dots, x_n),$
where
$$x_{j+1}=\cos(j\alpha_\tau)e_1 + \sin(j\alpha_\tau)e_2,\quad {\rm and}\quad
\alpha_\tau = \pi(n-1-2p(1-\tau))/n,$$
for $j=0, 1, 2, \dots, n-2$, while
$$x_n = \cos((n-1)\alpha_\tau)e_1 + \sin((n-1)\alpha_\tau)e_2 +
\sin(\pi\tau)e_3.$$
It is clear that $H_0=I_p$ and $H_1=I_0$.
Since $u|_{V_0}\ne 0$  we obtain $I_0^\ast(u)\ne 0$ as a consequence of
Proposition 10.3 in \cite{Br1}. This Proposition describes the cohomology
of Stiefel
manifolds with $\Z_2$ coefficients; it implies that $I_0^\ast:
H^\ast(V_p;\Z_2)\to H^\ast(W;\Z_2)$
is a monomorphism.
Hence, it follows $I_p^\ast(u) =I_0^\ast(u) \ne 0$ and thus $u|_{V_p}\ne 0$.

We want to show that for even $i$ and $j$  with $i+j=2p+2$
 the following relation holds in the ring $H^\ast(F'_{p+1};\Z_2)$:
\begin{eqnarray}
\sigma_i \sigma_{j} = \displaystyle{2p+2 \choose i}\sigma_{2p+2}+
\epsilon_{ij}\sigma_{2p}ue^{m-2}, \label{prod}
\end{eqnarray}
where $\epsilon_{ij}\in \Z_2$.
Note that for $r= 2(p+1)(m-1)$ the homomorphism
$$\Phi: H^{r}(F'_{p+1};\Z_2) \to H^{r}(F_{p+1};\Z_2)$$
has one-dimensional image and one-dimensional kernel.
For any even $r=2, 4, \dots, 2p+2$ the image of the Thom class $\sigma_{r}$
under $\Phi$
is nonzero and so it equals
the class $\Phi(\sigma_r)\in H^\ast(F_{p+1};\Z_2)$, which was denoted in
Theorem \ref{thm3}
by $\sigma_r$. From Theorem \ref{thm3} we know that
$$\Phi(\sigma_i)\Phi( \sigma_{j}) = \displaystyle{2p+2 \choose
i}\Phi(\sigma_{2p+2}).$$
This proves  (\ref{prod}) since $\sigma_{2p}ue^{m-2}$ belongs to the
kernel of $\Phi$.

The rest of properties in hypothesis $\mathcal F_{p+1}$ are now obvious.
Thus $\mathcal F_p$ implies $\mathcal F_{p+1}$, and
the proof of Theorem \ref{thm61} is complete. \proofend

\begin{remark}\label{rmk9} {\rm The above proof shows that the function
$L_X: G(X, n)\to \R$
is  perfect with respect to the field $\Z_2$ in the $D_n$-equivariant sense
as well, compare Proposition \ref{perf}.}\end{remark}

We conclude the paper with a proof of Theorem \ref{thm1} formulated in
Introduction.

{\bf Proof of Theorem \ref{thm1}.} Statement (B) follows from Proposition
\ref{thm7}
and Theorem \ref{thm61} that states that the sum of Betti numbers of the space
$G(S^m;n)/D_n$ is $m(n-1)$.

To prove statement (A) of Theorem \ref{thm1} one needs to
find the cup-length of $H^\ast(G(S^m,n)/D_n;\Z_2)$. The argument is similar
to the one in
Corollaries \ref{cor1} and \ref{cor2}. Namely, the length of the longest
nontrivial product
of classes $\sigma_{2i}$ from Theorem \ref{thm61}
equals the maximal number of 1's in the binary expansion, that {\it even}
integers,
not exceeding $n-3$, may have. This number is equal to the maximal number
of 1's in the
binary expansion, that integers, not exceeding $(n-3)/2$, may have, that
is, to
$$[\log_2((n-3)/2+1)] = [\log_2(n-1)]-1.$$
We also have the classes $u$ and $e^{m-1}$ at our disposal; therefore the
cup-length equals
$[\log_2(n-1)] + m -1$. Statement (A) follows now from Proposition
\ref{thm7}. \proofend

\end{document}